\documentclass[11pt]{amsart}

\DeclareMathOperator{\sech}{sech}
\DeclareMathOperator{\relu}{relu}
\DeclareMathOperator{\sigmoid}{sigmoid}
\DeclareMathOperator{\elu}{elu}

\usepackage{lipsum}
\usepackage{amsfonts,bm,amssymb,amsmath}
\usepackage{graphicx}
\usepackage{epstopdf}
\usepackage[caption=false]{subfig}
\usepackage{pgfplots}
\usepackage{algorithmic}
\usepackage{amsopn}
\usepackage{amsrefs}
\usepackage{float}
\usepackage{booktabs}

\usepackage[normalem]{ulem}
\usepackage{geometry}

\usepackage{hyperref}
\usepackage{cleveref}

\usepackage{color}
\usepackage{diagbox}

\theoremstyle{definition}

\newtheorem{example}{Example}[section]

\hypersetup{
	colorlinks,
	linkcolor={red!50!black},
	citecolor={blue!50!black},
	urlcolor={blue!80!black}
}

\numberwithin{equation}{section}

\definecolor{newcolor1}{rgb}{.8,.349,.1}
\colorlet{bblue}{blue!50!black}

\definecolor{newcolor1}{rgb}{.8,.349,.1}
\definecolor{airforceblue}{rgb}{0.36, 0.54, 0.66}
\definecolor{amethyst}{rgb}{0.6, 0.4, 0.8}
\definecolor{applegreen}{rgb}{0.55, 0.71, 0.0}
\definecolor{carrotorange}{rgb}{0.93, 0.57, 0.13}
\definecolor{ceruleanblue}{rgb}{0.16, 0.32, 0.75}
\definecolor{cornflowerblue}{rgb}{0.39, 0.58, 0.93}
\definecolor{darkpastelpurple}{rgb}{0.59, 0.44, 0.84}

\renewcommand{\i}{\ensuremath{\text{\normalfont I}}}
\newcommand{\ii}{\ensuremath{\text{\normalfont I\!I}}}
\newcommand{\iii}{\ensuremath{\text{\normalfont I\!I\!I}}}
\newcommand{\iv}{\ensuremath{\text{\normalfont I\!V}}}
\crefformat{equation}{(#2#1#3)}
\crefmultiformat{equation}{(#2#1#3)}{ and~(#2#1#3)}{, (#2#1#3)}{ and~(#2#1#3)}

\crefformat{figure}{Figure~#2#1#3}
\crefmultiformat{figure}{Figures~ #2#1#3}{ and~#2#1#3}{, (#2#1#3)}{ and~(#2#1#3)}
\crefformat{table}{Table~#2#1#3}

\def\Nm{\mbox{$\mathcal{N}$}}
\def\Wm{\mbox{$\mathcal{W}$}}
\def\Hm{\mbox{$\mathcal{H}$}}
\def\Lm{\mbox{$\mathcal{L}$}}

\def\k{\mbox{\boldmath $k$}}

\def\x{\mbox{\boldmath $x$}}

\def\0{\mbox{\boldmath $0$}}

\newcommand{\red}{\color{black}}

\begin{document}

\title[A machine-learning method for wave equations]{A machine-learning method for time-dependent wave equations over unbounded domains}

%\thanks{Support information for the second author.}

%    Address of record for the research reported here
%\address{Department of Mathematics, Louisiana State University, Baton
%Rouge, Louisiana 70803}
%    Current address
%\curraddr{Department of Mathematics and Statistics,
%Case Western Reserve University, Cleveland, Ohio 43403}
%\email{xyz@math.university.edu}
%    \thanks will become a 1st page footnote.

%\thanks{The first author was supported in part by NSF Grant \#000000.}
\author{Changjian Xie}
\address{School of Mathematical Sciences\\ Soochow University\\ Suzhou\\ China.}
\email{20184007005@stu.suda.edu.cn}

\author{Jingrun Chen}
\address{School of Mathematical Sciences\\ Soochow University\\ Suzhou\\ China.}
\email{jingrunchen@suda.edu.cn (Corresponding author)}

\author{Xiantao Li} 
\address{Department of Mathematics\\ The Pennsylvania State University\\ University Park \\ PA 16802\\ USA.}
\email{xiantao.li@psu.edu (Corresponding author)}
%\thanks{Support information for the second author.}

%    General info

\date{\today}

%\dedicatory{This paper is dedicated to our advisors.}


\begin{abstract}
Time-dependent wave equations represent an important class of partial differential equations (PDE) for describing wave propagation phenomena, which are often formulated over unbounded domains. Given a compactly supported initial condition, classical numerical methods reduce such problems to bounded domains using artificial boundary condition (ABC). In this work, we present a machine-learning method to solve this type of equations as an alternative to ABCs. Specifically, the mapping from the initial conditions to the PDE solution is represented by a neural network, trained using wave packets that are parameterized by their {band width} and wave numbers. The accuracy is tested for both  the second-order wave equation and the Schr\"{o}dinger equation, including  the nonlinear Schr\"{o}dinger equation. We examine the accuracy from both  \textit{interpolations} and \textit{extrapolations}. For initial conditions lying in the training set, the learned map has good interpolation accuracy, due to the approximation property of deep neural networks. The learned map also exhibits some  good extrapolation accuracy. % when the wave hits the boundary of the computational domain even though it does not ``see" the boundary in the training data set. 
We also demonstrate the effectiveness of the method for problems in irregular domains.
Overall, {the proposed} method provides an interesting alternative for finite-time simulation of wave propagation.
\end{abstract}

\maketitle

\section{Introduction}

Wave propagation is an ubiquitous phenomenon and for a long time, the associated properties have been a subject of interest in many disciplines \cite{whitham2011linear}. Aside from the well known acoustic waves,   the Schr\"{o}dinger equation that describes electronic waves, and  the elastodynamics  that embodies stress waves \cite{eringen2012electrodynamics} are also important examples. These models share the common ground that waves often propagate in an unbounded domain, even though they are triggered locally, e.g., by the presence of a wave source or an external forcing.
%The issue of unbounded domain poses challenges for numerical simulations.

One classical numerical approach to treat wave propagation in an unbounded domain is the absorbing boundary condition (ABC),
which confines the computation to a finite domain, and an ABC is then imposed on the boundary to minimize undesirable reflections \cites{arnold2003,Berenger1994,engquist1977absorbing,givoli2008computational,han2013artificial}. 
Rather than simply removing the exterior region, the ABC provides an efficient approach to mimic the influence from the surrounding environment. There are several different approaches to construct and implement ABCs, most of which involve the derivation and approximation of the Dirichlet-to-Neumann map.  There has been a large body of works on ABCs and interested readers may refer to the review articles \cites{antoine2007review,givoli2008computational} for details and references therein. %Among works on the construction of ABCs \cite{engquist1977absorbing, fevens1999absorbing, antoine2006artificial, xu2006absorbing, antoinenonlinear2008review,givoli2008computational, antoine2011absorbing, li2019absorbing, WuLi2020}, most efforts have been focused on one-dimensional problems with a few exceptions for multi-dimensional problems with flat boundary \cite{antoine2004numerical, han2004exact,szeftel2004design,xu2007adaptive, jiang2008efficient}. 
The integration of ABCs with finite difference or finite element methods, has also been extensively studied \cites{shibata1991absorbing, moxley2013generalized, li2018efficient}.

Recently, the rapid progress in deep learning has driven the development of solution techniques for PDEs under the framework of deep learning, especially in high-dimensional cases where deep neural networks (DNNs) are expected to overcome the curse of dimensionality; see \cite{E2020} for a review and \cites{lagaris1998artificial, E2017Dec,khoo2017solving, E2018Mar,han2018solving,sirignano2018dgm,raissi2019physics-informed,Beck2019Aug,Hutzenthaler2019Mar,cai2019phase, zang2020weak,Becker2020May, Lyu2020Jun} for specific examples. 
{One remarkable application of neural networks is the physics-informed neural networks (PINNs)  \cites{raissi2019physics-informed}, which has demonstrated its accuracy in solving both forward problems and inverse problems,  where model parameters are inferred from the observed data. PINNs have already been applied to a range of problems, including those in fluid dynamics \cites{Raissi2020HiddenFM,Jin2020NSFnetsF}, meta-material design \cites{Chen2020PINN}, biomedical engineering \cites{Yazdani2019SystemsBI}, uncertainty quantification \cites{YANG2021BPINN} and free boundary problems, besides the high dimensional PDEs and stochastic differential equations.} Typically, the loss function is defined over a {\it finite} domain in most methods, such as the deep Ritz method \cite{E2018Mar}, deep Galerkin method \cite{sirignano2018dgm}, physics-informed neural networks \cite{raissi2019physics-informed}, and deep mixed residual method \cite{Lyu2020Jun}. To the best of our knowledge, the only exceptions are the full history recursive multilevel Picard approximation method \cites{Beck2019Aug,Hutzenthaler2019Mar} and the deep backward stochastic differential equation method \cites{E2017Dec,han2018solving}, where the solution of the underlying PDE is approximated through the solution of a suitable stochastic optimization problem on an appropriate function space. Typical equations are (semilinear) parabolic PDEs.
These recent works have demonstrated the strong representability of DNNs for solving PDEs.

The current work aims to solve time-dependent wave equations on {\it unbounded} domains using deep learning. One natural approach is to build an artificial neural network (ANN) that takes an  ABC, e.g., the perfectly matched layer (PML) method, into account \cites{yao2018machine, yao2020enhanced}. { The basic idea behind the approach in \cites{yao2018machine} is as follows. Given the electromagnetic field at the current step, one can predict the field on the PML boundary at next time step. Then, one computes a field in a slightly larger domain, called the object domain, at the next step through output from PML, which subsequently becomes the new input data, by the finite-difference time-domain (FDTD) method. Furthermore, one can embed the network model into the FDTD method and replace the PML. The data groups are collected at the interface with conventional PML. The Long Short Term Memory (LSTM) network based on the PML model in \cites{yao2020enhanced} can achieve higher accuracy than an ANN that is based on the PML model, thanks to the sequence dependence feature of LSTM networks. Compared to the conventional PML approach,  the machine-learning methods in \cites{yao2018machine,yao2020enhanced} decrease the size of the boundary region and the complexity of the FDTD method, due to the introduction of a one-cell boundary layer. But the data generation involves prior PML computation.} This process involves the history of solutions at the boundary and may be rather complicated in  general. We propose a different machine-learning strategy to solve the {wave propagation over unbounded domain}. Given a compactly supported initial condition, we restrict the full problem to a solution mapping over a finite region. More specifically, the mapping from the initial condition, expressed as wave packets with  {band width} and wave numbers as parameters, to the PDE solution in the same compactly supported domain at later times, is represented by a fully connected neural network  {or residual neural network}. The parameters in the network are then trained using data that can be generated using a variety of methods, ranging from   analytical solutions,  numerically computed solutions, to  approximated solution from PINNs.% One key point of the training data set is that the collected data correspond to temporal instances that the wave does not hit the boundary, i.e., no ABC is needed to generate the training data in practice, although theoretically the ABC is required. Surprisingly, such a map has good \textit{interpolative} property and \textit{extrapolative} property.

On one hand, the mapping can generate accurate results in which the specific initial condition is not included in the training set, but can be interpolated by those in the training set. On the other hand, the method also allows extrapolations, e.g.,  when the wave packet arrives at the boundary of the finite region, even though the training set  only contains temporal instances prior to that event. Compared to existing works,  {the proposed} method can be easily implemented. The solutions represented by DNN also exhibit absorbing properties. But there is no need to determine the coefficients in ABCs, or to incorporate ABCs into finite difference or finite element methods.
{The proposed} method provides an alternative for finite-time simulation of wave propagation.

This paper is organized as follows. In \cref{sec:main method}, we describe the machine-learning method for two representative wave equations:  the second-order wave equation and the Schr\"{o}dinger equation. Numerous examples are provided to show the interpolative and extrapolative properties of the proposed method in \cref{sec:experiments}. Conclusions are drawn in \cref{sec:conclusions}. 

\section{Methodology}\label{sec:main method}

To elaborate the approach of constructing solution representations by a neural network, we consider, as specific examples, the time-dependent acoustic wave equation and the Schr\"odinger equation as examples,  due to the fact that they have been treated extensively in PDE analysis and numerical approximations. But we expect that the idea can be extended to other types of wave equations. We express these two models  as time-dependent PDEs  over the entire space $\mathbb{R}^d$:
\begin{enumerate}
	\item[(I) ] Time-dependent wave equation:
	\begin{align}\label{wave-eq}
	\begin{aligned}
	u_{tt}&=\Delta u, \quad \x \in \mathbb{R}^d,\;t>0,\\
	u(\x,0)&=u_0(\x),\quad u_t(\x,0)= v_0(\x).
	\end{aligned}
	\end{align}
	\item[(II) ] Time-dependent Schr\"odinger equation: 
	\begin{align}\label{Schrodinger-1}
	\begin{aligned}
	i \partial_t u(\x,t)&=- \Delta u(\x,t)+V(\x,t)u(\x,t)+f(|u(\x,t)|^2)u(\x,t),\quad \x \in \mathbb{R}^d,\quad t>0,\\
	u(\x,0)&=u_0(\x),\quad \x \in \mathbb{R}^d.
	\end{aligned}
	\end{align}
	
	There are important cases that deserve particular attention:
	\begin{enumerate}
		\item The linear Schr\"odinger equation: $f\equiv 0$. This describes the dynamics of a free electron.
		\item The cubic Schr\"odinger equation:
		\begin{align}\label{nonlinear}
		f(\rho)=\beta \rho, \quad \rho \in [0,\infty),
		\end{align} 
		where $\beta$ (positive for repulsive or defocusing interaction and negative for attractive or focusing interaction) is a given dimensionless constant describing the strength of the interaction.   It has been widely used to  model nonlinear wave interactions in a dispersive medium. 
	\end{enumerate}
\end{enumerate}

We have expressed these models in their non-dimensionalized forms. For example, the wave speed in \eqref{wave-eq} and the Planck constant in \eqref{Schrodinger-1} have been both set to unity.  In addition, we set   $\beta=-1$ in \cref{nonlinear}.

\medskip
We make the important assumption that the initial condition and the potential are compactly supported in a finite domain, denoted by $\Omega \subset \mathbb{R}^d$, that is,
\[\text{supp}(u_0),~ \text{supp}(v_0),~  \text{supp}(V)  \subset \Omega.\] 

Our aim is to determine the solution $u(\cdot,t)$ in the same domain $\Omega$ at later times.

\subsection{The training procedure}

In this section, we describe how the solution is trained using neural networks.
One key step in a machine learning procedure is the preparation of a dataset, which will subsequently be fed into the machine learning model.
We first prepare a dataset, consisting of the initial condition and the corresponding solutions at later times. We will denote initial data by $U_0=(u_0,\partial_t u_0)$ for time-dependent wave equation \cref{wave-eq} and { $U_0=(p_0, q_0)$ for the time-dependent Schr\"odinger equation \cref{Schrodinger-1} with $p$ and $q$ being the real  and imaginary parts of the wave function, respectively.} In principle, the mapping from $U_0$ to the solution at a later time can be expressed as an operator $\mathcal{S}$, 
\begin{equation}
u(\x,t)|_\Omega = \mathcal{S} U_0. 
\end{equation}
For example, in the linear case, this can be written as an integral operator using the Green's function \cite{Evans}. But such an expression is of limited value in practice since the direct evaluation is rather expensive. Here we represent such a mapping  using a neural network and determine the parameters through training.

In the training step, we consider three cases, as motivated by the terminology in control systems, 
\begin{enumerate}
	\item[(I) ] Single-input single-output (SISO) datasets
	\begin{equation}\label{eq: dataset-single}
	%U_T=\{u_0^\ell,\partial_t u_0^\ell,u_T^\ell\}_{\ell=1}^N, \textrm{ for (M1), and } U_T=\{u_0^\ell,u_T^\ell\}_{\ell=1}^N,   \textrm{ for (M2)}.
	\left\{U_0^\ell,U_T^\ell\right\}_{\ell=1}^N.
	\end{equation}
	\item[(II) ] Single-input multiple-output (SIMO) datasets
	\begin{equation}\label{eq: dataset-multi}
	% 		 U_{\{t_i\}_{i=1}^p}=\{u_0^\ell,\partial_t u_0^\ell,u_{\{t_i\}_{i=1}^p}^\ell\}_{\ell=1}^N, \textrm{ for (M1), and } U_{\{t_i\}_{i=1}^p}=\{u_0^\ell,u_{\{t_i\}_{i=1}^p}^\ell\}_{\ell=1}^N, \textrm{ for (M2)}.   
	\left\{U_0^{\ell},U_{\{t_i\}_{i=1}^p}^\ell\right\}_{\ell=1}^N
	\end{equation}
	\item[(III) ] Exogenous-input multiple-output (XIMO) datasets
	\begin{equation}\label{eq: dataset-ex}
	\left\{V_{\{t_i\}_{i=1}^p}^\ell,U_{\{t_i\}_{i=1}^p}^\ell\right\}_{\ell=1}^N.
	\end{equation}
\end{enumerate}
%or extended dataset, denoted by
%\begin{equation}\label{eq: dataset}
%U_{\{t_i\}_{i=1}^p}=\{\psi_0^\ell,\psi_{\{t_i\}_{i=1}^p}^\ell\}_{\ell=1}^N,    
%\end{equation}
Here the integer $N$ denotes the number of training samples, and $p$ refers to the time instances where the solutions are observed. 
The {\it input} simply refers to the initial conditions and the output involves the resulting solutions at a later time (or at multiple time instances).  Namely,  $U_T:=\left\{ u(\cdot, T)\right\}. $ These solutions will be collected at grid points that lie in the domain of interest $\Omega$. For simplicity, we also work with $U_T^\ell$ in the same domain. But in practice,  one can also choose $U_T$ in a different domain.
In the case of XIMO, one may consider the Schr\"odinger equation, with the initial condition fixed at ground state. The dynamics  is then entirely driven by the external potential.

% The initial condition will be considered here as Gaussian wave packet type with compact support. 
%  defined by
% \[\text{supp}(g):=\overline{\{x\in R| g(x)\neq 0\}} \subset \Omega\] in which $\Omega$ is bounded domain. 

%A natural example $\psi_0$ may be  Gaussian wave packets  and $\psi_T$ being analytical   or numerical solution at final time $t=T$. 

\begin{figure}[htbp]
	\begin{center}
		\begin{tikzpicture}
		%		\draw [draw=none,rounded corners,fill=yellow!20](-8.2,0.2) rectangle (6.2,-17);
		\draw [draw=none,rounded corners,fill=blue!20](-8,-4.5) rectangle (-2,-6);
		\draw [draw=none,rounded corners,fill=blue!20](0.3,-4.5) rectangle (5.9,-6);
		\node [draw=none,rounded corners,fill=orange!50,text width=2.3cm](a)at(-5,-1){$U_0(\x,\sigma^2_1,k_1)$};
		\node (b)at(-5,-1.5){$\cdots$};
		\node [draw=none,rounded corners,fill=orange!50,text width=2.3cm](a)at(-5,-2){$U_0(\x,\sigma^2_1,k_J)$};
		\node [draw=none,rounded corners,fill=orange!50,text width=2.3cm](a)at(-5,-3){$U_0(\x,\sigma^2_2,k_1)$};
		\node (b)at(-5,-3.5){$\cdots$};
		\node [draw=none,rounded corners,fill=orange!50,text width=2.3cm](a)at(-5,-4){$U_0(\x,\sigma^2_K,k_J)$};
		\node (h1)at(-5,-4.8){Input};
		\node (h2)at(-5,-5.2){Initial condition of};
		\node (h3)at(-5,-5.6){Gaussian wave packets};
		
		\node [draw=none,rounded corners,fill=orange!50,text width=2.3cm](a1)at(3,-1){$U_T(\x,\sigma^2_1,k_1) $};
		\node (b)at(3,-1.5){$\cdots$};
		\node [draw=none,rounded corners,fill=orange!50,text width=2.3cm](a1)at(3,-2){$U_T(\x,\sigma^2_1,k_J) $};
		\node [draw=none,rounded corners,fill=orange!50,text width=2.3cm](a1)at(3,-3){$U_T(\x,\sigma^2_2,k_1) $};
		\node (b)at(3,-3.5){$\cdots$};
		\node [draw=none,rounded corners,fill=orange!50,text width=2.3cm](a1)at(3,-4){$U_T(\x,\sigma^2_K,k_J) $};
		\node (hh1)at(3,-4.8){Output};
		\node (hh2)at(3,-5.2){from exact, numerical solution};
		\node (hh3)at(3,-5.6){or neural network prediction};
		
		\draw [->,line width=1.4pt,color=applegreen](-3.7,-1) -- (1.8,-1);
		\draw [->,line width=1.4pt,color=applegreen](-3.7,-2) -- (1.8,-2);
		\draw [->,line width=1.4pt,color=applegreen](-3.7,-3) -- (1.8,-3);
		\draw [->,line width=1.4pt,color=applegreen](-3.7,-4) -- (1.8,-4);
		
		\draw [->,line width=1.4pt,color=applegreen](-3.7,-1) -- (1.8,-2);
		\draw [->,line width=1.4pt,color=applegreen](-3.7,-1) -- (1.8,-3);
		\draw [->,line width=1.4pt,color=applegreen](-3.7,-1) -- (1.8,-4);
		
		\draw [->,line width=1.4pt,color=applegreen](-3.7,-2) -- (1.8,-1);
		\draw [->,line width=1.4pt,color=applegreen](-3.7,-2) -- (1.8,-3);
		\draw [->,line width=1.4pt,color=applegreen](-3.7,-2) -- (1.8,-4);
		
		\draw [->,line width=1.4pt,color=applegreen](-3.7,-3) -- (1.8,-1);
		\draw [->,line width=1.4pt,color=applegreen](-3.7,-3) -- (1.8,-2);
		\draw [->,line width=1.4pt,color=applegreen](-3.7,-3) -- (1.8,-4);
		
		\draw [->,line width=1.4pt,color=applegreen](-3.7,-4) -- (1.8,-1);
		\draw [->,line width=1.4pt,color=applegreen](-3.7,-4) -- (1.8,-2);
		\draw [->,line width=1.4pt,color=applegreen](-3.7,-4) -- (1.8,-3);
		
		\draw [->,line width=1.4pt,color=red!10!green!20!blue!60](-2,-5) -- (0.3,-5);
		\node (hh4)at(-1,-4.5){$\Nm_D^M\textrm{ or } \tilde{\Nm}_D^M$};
		\node (hh5)at(-1,-5.5){Training};

		\draw [->,line width=1.4pt,color=red!10!green!20!blue!60](3,-6) -- (3,-6.7) -- (-4,-6.7);
		\node (hh5)at(-5.3,-6.7){Predict $U_T(\x)$};
		\node (hh4)at(-1,-6.4){Given $U_0(\x)$};
		\node (hh5)at(-1,-7){Testing};
		
		%		\node (cap1)at(-1,-7.5){\text{(A1) Learning framwork $\Nm_D^M(u_0,\partial_tu_0,\Wm)$ or $\tilde{\Nm}_D^M(u_0,\partial_tu_0,\Wm)$ }};
		%		\node (cap1)at(-1,-8){with single-snapshot output for time-dependent wave equation.};
		
		%			\node (cap1)at(-1,-14.5){\text{(B) Learning framwork $\Nm_D^M(\psi_0,\Wm)$ for time-dependent Schr\"odinger equation.}};
		
		%		\draw [->,line width=1.4pt,color=red!10!green!20!blue!60](3,-14.5) -- (3,-15.2) -- (-4,-15.2);
		%		\node (hh5)at(-5.3,-15.2){Predict $u(\x,T)$};
		%		\node (hh4)at(-1,-14.9){Given $u_0(x)$};
		%		\node (hh5)at(-1,-15.5){Testing};
		
		%		\node (cap1)at(-1,-16){\text{(B1) Learning framwork $\Nm_D^M(u_0,\Wm)$ or $\tilde{\Nm}_D^M(u_0,\Wm)$}};
		%		\node (cap1)at(-1,-16.5){with single-snapshot output for time-dependent Schr\"odinger equation.};
		\end{tikzpicture} 
		
		\caption{Schematic of learning mapping $\Nm_D^M$ (FCNNs) or $\tilde{\Nm}_D^M$ (ResNets) from initial condition of Gaussian wave packets with compact support to the output data from exact, or numerical solution, or neural network prediction for the time-dependent wave equations.}\label{fig-mode1}
	\end{center}
\end{figure}
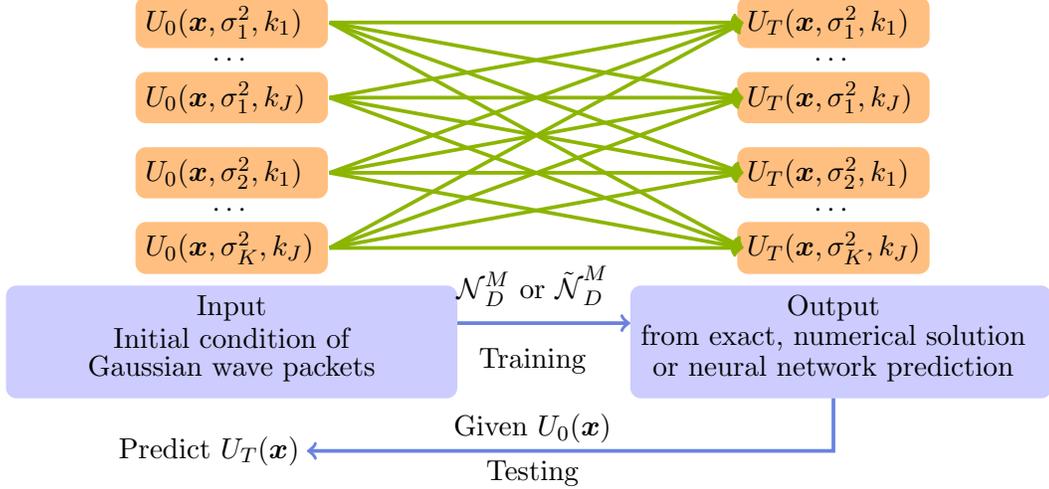

We will discuss the construction of datasets in the next section in more details. In particular, properties of wave propagations, e.g., wave length and dispersion relations, are built into the training set.   % Specifically, the initial condition and numerical solution can be parametrized by $\psi_0(x,\sigma^2,k)$ and $\psi_T(x,\sigma^2,k)$.
{ Our experience suggests that properly rescaling  the input data can  improve the convergence. Specifically, before it is fed into the network, the initial condition will go through the following transformation.
	\begin{align*}
	{U}_0 \rightarrow\frac{2(\lambda U_0-\inf(U_0))}{\sup(U_0)-\inf(U_0)}-\inf(U_0),
	\end{align*}
	where $\inf(U_0)$ and $\sup(U_0)$ denote the infimum and supremum of $U_0$, respectively. Note that after the rescaling, the input data ${U}_0$ take values in  $[-1,1]$.}
Since $U_T$ is fully determined by $U_0$, we approximate the mappings from $U_0$ to $U_T$ using a neural network, denoted by $\Nm_{D}^M(U_0,\Wm)$, i.e.,
\begin{equation}\label{eq: nnet0}
U_T \approx \Nm_{D}^M(U_0,\Wm).
\end{equation}
{ The neural network underlying the mapping \eqref{eq: nnet0} is illustrated in \cref{fig-mode1}.}

The function $\Nm_{D}^M$ is determined by a network consisting of $D$ layers with width $M$, and the associated parameters are denoted by $\Wm$. %With $\Nm$ being a functional which indicated the functional learning, 
For a fully connected neural network (FCNN), the mapping \cref{eq: nnet0} from input to output is explicitly given by
\begin{align}\label{eq: nnet1}
\begin{aligned}
\Nm_D^M(U_0,\Wm)&=W_D^T\Hm_{D-1}^M(U_0,\tilde{\Wm})+b_D,
\end{aligned}
\end{align}
where,
\begin{align*}
\begin{aligned}
\Hm_{D-1}^M(U_0,\tilde{\Wm})&=\phi(W_{D-1}^T\cdots \phi(W_2^T\phi(W_1^T U_0+b_1)+b_2)\cdots+b_{D-1}),
\end{aligned}
\end{align*}
with $\phi$ being the activation function and $\{W_j,b_j\}_{j=1}^D$ being the parameters specified by the network.

The residual neural network (ResNet) structure \cite{He2015} will also be considered in our numerical studies. In this case, the mapping can be expressed with the following steps, 
\begin{align*}%\label{eq: nnet}
% 	\begin{cases}
% 	y=W_1^T[u_0,\partial_t u_0]+b_1\\
% 	y=y+\Hm_{D_1}^M(y)\\
% 	y=y+\Hm_{D_2}^M(y)\\
% 	\cdots\quad  \cdots \quad \cdots\\
% 	y=y+\Hm_{D_s}^M(y)\\
% 	\tilde{\Nm}_D^M (u_0,\partial_t u_0,\Wm)= W_D^Ty+b_D,
% 	\end{cases}
% 	\textrm{ for (M1), and }
\begin{cases}
y_1=W_1^TU_0+b_1,\\
y_2=y_1+\Hm_{D_1}^M(y_1),\\
y_3=y_2+\Hm_{D_2}^M(y_2),\\
\cdots\quad  \cdots \quad \cdots\\
y=y_s+\Hm_{D_s}^M(y_s),\\
\tilde{\Nm}_D^M (U_0,\Wm)= W_D^Ty+b_D,
\end{cases}
\end{align*}
%\begin{align*}%\label{eq: nnet}
%\begin{aligned}
%\Nm_D^M(u_0,u_1,\Wm)&=W_D^T h(W_{D-1}^T\cdots h(W_2^Th(W_1^T[u_0,u_1]+b_1)+b_2)\cdots+b_{D-1})+b_D,\textrm{ for (M1)},\\
%\Nm_D^M(\psi_0,\Wm)&=W_D^T h(W_{D-1}^T\cdots h(W_2^Th(W_1^T\psi_0+b_1)+b_2)\cdots+b_{D-1})+b_D,\textrm{ for (M2)}
%\end{aligned}
%\end{align*}
where $s$ is the number of residual blocks with a skip connection. 

%The neural network underlying the mapping \eqref{eq: nnet0} is illustrated in \cref{fig-mode1}.% Also shown in the diagram is the case where solutions at multiple instances $\{t_i\}_{i=1}^p$ are included as the output.  

\medskip

The next step is to formulate the problem as a supervised learning problem by means of minimizing the population risk (expected risk), elabrated in \cites{bottou2018optimization} by  
% {\color{red} Can we cite a standard reference?}
\begin{align*}
\begin{aligned}
\min_{\Wm} \mathbb{E}_{(u_0,u_T)\sim \mu}\left[\|\Nm_D^M(u_0,\Wm)-u_T\|^2\right],\textrm{ or } \min_{\Wm} \mathbb{E}_t \mathbb{E}_{(u_0,u(t))\sim \mu}\left[\|\Nm_D^M(u_0,\Wm)-u(t)\|^2\right],
\end{aligned}
\end{align*}
with $\mu$ being a probability distribution, which in practice, can be discretized by mean squared error as empirical loss for the training samples. For example, for a  {SISO} dataset, this leads to a cost function, 
{\small
	\begin{align}\label{eq: training}
	%	\begin{aligned}
	\Lm(U_0,\Wm)=	\frac{1}{N}\sum_{\ell=1}^N\Big[\Nm_D^M(U_0^\ell,\Wm)-u_T^\ell\Big]^2.
	%	\end{aligned}
	\end{align}}

Similarly, for a  {SIMO} dataset, we can define the loss function as follows,
{\small
	\begin{align}\label{eq: training-simo}
	\Lm(U_0,\Wm)=\frac{1}{N}\sum_{i=1}^p\sum_{\ell=1}^N\Big[\Nm_D^M(U_0^\ell,\Wm)-u_{\{t_i\}_{i=1}^p}^\ell\Big]^2.
	\end{align}}

%\begin{equation}\label{eq: training}
%\min_{\Wm} \frac{1}{N} \sum_{\ell=1}^N\|\Nm_D^M(\psi_0^\ell,\Wm)-\psi_T^\ell\|^2, \textrm{ or } \sum_{\ell=1}^N\|\Nm_D^M(\psi_0^\ell,\Wm)-\psi_{\{t_i\}_{i=1}^p}^\ell\|^2.
%\end{equation}
% {\red we need to explain this part more
% However, we extremely focus on the population risk for the performance of the unseen sample as
% \begin{align*}
% 	\begin{aligned}
% % 	\min_{\Wm} \mathbb{E}_{(u_0,\partial_t u_0,u_T)\sim \mu}\left[\|\Nm_D^M(u_0,\partial_t u_0,\Wm)-u_T\|^2\right],\textrm{ or } \min_{\Wm} \mathbb{E}_t \mathbb{E}_{(u_0,\partial_t u_0,u(t))\sim \mu}\left[\|\Nm_D^M(u_0,\partial_t u_0,\Wm)-u(t)\|^2\right], \\
% 	\min_{\Wm} \mathbb{E}_{(u_0,u_T)\sim \mu}\left[\|\Nm_D^M(u_0,\Wm)-u_T\|^2\right],\textrm{ or } \min_{\Wm} \mathbb{E}_t \mathbb{E}_{(u_0,u(t))\sim \mu}\left[\|\Nm_D^M(u_0,\Wm)-u(t)\|^2\right],
% 	\end{aligned}
% \end{align*}
% with $\mu$ being a probability distribution. 
% }

%\subsection{Notations}
\medskip
Notations for the parameters of our model and algorithm are summarized in \cref{tab:notation}.  
\begin{table}[htbp]
	\centering
	\begin{tabular}{|c|c|}
		\hline
		$d$ & the dimension of the problem  \\
		\hline
		$\mathcal{K}$ & the set of wave numbers of wave packet \\
		\hline
		$\Sigma$ & the set of width of wave packet \\
		\hline
		$D$ & number of layers \\
		\hline
		$M$ & number of neurons of each hidden layer \\
		\hline
		$N$ & number of initial conditions (training samples) \\
		\hline
		$m_1$ & number of column of input matrix \\
		\hline
		$m_2$ & number of column of output matrix \\
		\hline
		$p$ & number of time instances \\
		\hline
		$(\cdot)_{F}=\Nm_D^M$ & representation of fully connected neural networks \\
		\hline
		$(\cdot)_{R}=\tilde{\Nm}_D^M$ & representation of residual neural networks \\
		\hline
		$\beta$ & constant describing the strength of interaction \\
		\hline
		$\lambda$& rescaling parameter of input data\\
		\hline
		$s$ & number of residual blocks \\
		\hline
		$E$ & exponential spacing with $\Sigma=\{h,2h,2^2h,2^3h,2^4h,2^5h\}$\\
		\hline
		$L$ & linear spacing with $\Sigma=\{0.8,0.9,1,1.1,1.2,1.3\}$\\
		\hline 
		$\delta t$ & time step \\
		\hline
		$p$ & the real part of wave function\\
		\hline
		$q$ & the imaginary part of wave function \\
		\hline
		$\mathcal{R}$& the $L^2$ relative error between the reference and DNN solution \\
		\hline
	\end{tabular}
	\caption{Notations for the parameters in the model  \eqref{eq: nnet1} and algorithm.}
	\label{tab:notation}
\end{table}

In addition to the structure of the network, another factor that may play a significant role in the approximation \cref{eq: nnet0} is the choice of the activation function.  In this paper, we first pick the FCNNs with the widely used \textit{relu} activation function. Then we implement a number of other nonlinear activation functions to test their accuracy, including:
\begin{align*}
&\begin{aligned}
\relu(x)=x_{+}=\begin{cases}
x, \quad x>0,\\
0, \textrm{ otherwise},
\end{cases}
\end{aligned}
\begin{aligned}
\tanh(x)=\frac{\exp(x)-\exp(-x)}{\exp(x)+\exp(-x)}
\end{aligned}\\
&\begin{aligned}
\sigmoid(x)=\frac{1}{1+\exp(-x)},
\end{aligned}\quad
\begin{aligned}
\elu(x)=\begin{cases}
x, \quad x>0,\\
\alpha(\exp(x)-1), \textrm{ otherwise}.
\end{cases}
\end{aligned}
\end{align*}

Our remaining task is to {\bf (i)} collect a suitable set of initial conditions and the corresponding solutions at later times, so that the  representation \cref{eq: nnet0} can be trained; {\bf (ii)} test the approximation   \cref{eq: nnet0} against analytical or numerical solutions.

\subsection{The training set for time-dependent wave equation}
Typical analysis of wave propagations starts with their dispersion properties using
Fourier transform \cite{whitham2011linear}.
For the time-dependent wave equation \cref{wave-eq}, 
the dispersion relation is given by $\omega(\k)=|\k|$ with $\k$ being the wave number. Often observed in practice are wave packets that are confined by an envelop and travel as a unit. Here we use wave packets  to form the training set. More specifically, we consider those wave packets with a Gaussian envelop, which can be derived, e.g.,  by using Fourier transform. For instance, for the acoustic wave equation \cref{wave-eq}, from the initial conditions 
\begin{align}\label{wave-initial-train-u}
u(\x,0)=\exp\left(-\frac{|\x|^2}{2\sigma^2}\right)\cos(\k\cdot \x),\quad \x \in \mathbb{R}^d, \quad \k \in \mathcal{K},\;\sigma^2 \in \Sigma,
\end{align} 
one obtains, 
\begin{align}
u(\x,t)=\exp\left(-\frac{\sum_{i=1}^d(x_i-t)^2}{2\sigma^2}\right)\cos(\k\cdot \x-|\k|t),\quad \x \in \mathbb{R}^d.
\end{align}

This implies that the initial velocity is given by,
\begin{align}\label{wave-initial-train-ut}
u_t(\x,0)=\exp\left(-\frac{|\x|^2}{2\sigma^2}\right)\left[\frac{\sum_{i=1}^d x_i}{\sigma^2} \cos(\k\cdot \x)+|\k|\sin(\k\cdot \x)\right],\quad \x \in \mathbb{R}^d.
\end{align}

%This wave packet solutions move in a direction by wave number $\k$ starting from the position at the origin in 1D case and move around the $xy$-plane that depends on $\k=(k_1,k_2)$ in 2D case. Generally, this type of solutions goes to anywhere in the hyperplane associated with the wave number $\k$ in high dimensions. 

\subsection{The training set for time-dependent Schr\"odinger equation}

\subsubsection{The linear case}\label{Dataset:linear}

%For any dimension $d$, we take $\varepsilon=1/2$, $V(x)=0$ and $f(\rho)=0$ in \cref{Schrodinger-1} 

For the linear Schr\"odinger equation,
\begin{align}\label{S-2}
i u_t = - \Delta u, \quad \x\in \mathbb{R}^d,\quad t>0,
\end{align}
the dispersion relation is given by $\omega(\k)=1/2 |\k|^2$.  To form the training set, we first pick the initial conditions from a family of wave packets,
\begin{align}\label{eq-p1}
u_0(\x)=\exp\left[-(\x/\sqrt{2})^2/\sigma^2+i\k\cdot (\x/\sqrt{2})\right],\quad \k\in \mathcal{K}, \; \sigma^2 \in \Sigma,
\end{align}
representing a Gaussian wave packet centered at the origin with wave number $\k$. The width parameter of the Gaussian envelope will be drawn from a pre-selected set: \[\sigma \in \Sigma, \; \Sigma\subset \Sigma_0 (:=\mathbb{R}_+), \quad \k \in \mathcal{K},\; \mathcal{K}\subset \mathcal{K}_0 (:=\mathbb{R}^d).\] 

%It is worthwhile to emphasize that $\mathcal{K}$ should be contain values that represent multiple directions in higher dimensions. For example, in  2D, $\k=(k_1,k_2)$, more samples for $k_1$ and $k_2$  usually lead to a more  robust datasets to be trained in the $xy$-plane. 
Theoretically, we can go through the whole spaces $\mathcal{K}_0$ and $\Sigma_0$. For practical purposes, we take some representative elements from finite subsets $\mathcal{K}$ and $\Sigma$. We will discuss more details about the selection of $\Sigma$ and $\mathcal{K}$ in the next section, and demonstrate how they impact the accuracy.

For each $\sigma$ and $\k$, the exact solution to \cref{S-2} can be constructed directly, 
\begin{align}\label{eq-1D}
u(x,t)=\frac{1}{\sqrt{1+2i t}}\exp\left(-\frac{1}{1+4t^2} \frac{(\frac{x}{\sqrt{2}}-kt)^2}{\sigma^2}\right) \exp\left(i\frac{1}{1+4t^2}[(k+\frac{2tx}{\sqrt{2}})\frac{x}{\sqrt{2}}-\frac12k^2t]\right), 
%\textrm{ in 1D},
\end{align}
for one-dimensional problems ($d=1$).

For two-dimensional problems ($d=2$), we have,
\begin{align}
u(x_1,x_2,t)=\left(\frac{i}{i-2t} \right) \exp\left[\frac{-i\left(\frac{x_1^2+x_2^2}{2\sigma^2}\right)-\frac{1}{\sqrt{2}}(k_1x_1+k_2x_2)+\frac12 (k_1^2+k_2^2)t}{i-2t}\right]. %\quad \textrm{ in 2D}.
\end{align}

These formulas can be generalized to arbitrary dimensions $d$, 
\begin{align}
u(\x,t)=\left(\frac{i}{i-2t}\right)^{\frac{d}{2}}\exp\left[\frac{-i\left(\frac{|\x|^2}{2\sigma^2}\right)-\frac{1}{\sqrt{2}}\k\cdot \x+\frac12|\k|^2t}{i-2t}\right].
%, \quad \textrm{ in dD}.
\end{align}

This family of solutions will constitute the datasets defined in \eqref{eq: dataset-single} and \eqref{eq: dataset-multi}, which will be used in the training step \eqref{eq: training}.% in this case can be specified by $\{\psi_0(\x),\psi(\x,t=T)\}$ in which $\sigma^2$ and $\k$ can be selected into the training process.

\subsubsection{The nonlinear case}\label{Dataset:nonlinear}

Since analytical solutions are difficult to obtain for the nonlinear PDE \cref{Schrodinger-1}, we generate the solutions $U_T^\ell$ using numerical methods. Here we use the finite difference scheme with uniform grid size, together with an operator-splitting scheme in time \cite{bao2002time}. More specifically, we use the Strang splitting, which at each  time step, involves the following operations:
\begin{enumerate}
	\item[(a) ] Solve $$i u_t+\Delta u=0,$$ 
	for half of the time step: $\delta t/2$.  Due to the linearity, this can be done exactly using the Fourier transform to diagonalize the Laplacian term. % $\psi(\x,\Delta t/2)$, with $\psi(\x,0)$ as initial data;
	\item[(b) ] Using the solution from the previous step, solve \[i u_t+|u|^2u=0,\] for one step.
	Using the fact that $\frac{d}{dt}|u|^2=0,$ this equation can also be solved exactly. This can also be extended to include an external scalar potential, that is \[i u_t+(|u|^2+V(x,t))u=0,\] for one step.
	\item[(c) ] Solve $iu_t+\Delta u=0$ again for another half step.
\end{enumerate}

% {\red what do we do for the case with time dependent potential?}

The symmetric operator splitting is known to have second order accuracy in time. In principle, one can also use higher order methods \cite{yoshida1990construction}, but the current numerical method is already adequate to test the neural network approximation. We  also pick initial conditions from \eqref{eq-p1}. The solutions at time $T$, together with the initial conditions \eqref{eq-p1}, will form the data set. 

% {\red How about the convergence rate of the SGD depends on $\mathcal{K}$ and $\Sigma$?}

%\begin{remark}
Another interesting approach to build the training set  is to design an embedding neural network to obtain the solutions to PDEs; see \cites{E2018Mar,Lyu2020Jun} for examples. {We will illustrate this approach using the nonlinear Schr\"odinger equation as an example. We follow the discrete-time PINNs \cite{raissi2019physics-informed} method, combined with the Crank-Nicolson finite difference scheme (CNFD) \cite{ANTOINE20132621} in time for \cref{Schrodinger-1} with $f(\rho)=\beta \rho^\mu$ and $\beta=-1$. Specifically, the discrete-time model is given by
	\begin{align}\label{eq-CNFD}
	i\frac{u^{n+1}-u^n}{\delta t}=-\frac12 \Delta\left(u^{n+1}+u^n\right)-{\frac14}\left[|u^{n+1}|^{2\mu}+|u^n|^{2\mu}\right](u^{n+1}+u^n),
	\end{align}
	which can be rewritten, using the representation $u=p+iq$, as follows,
	\begin{align}\label{eq-I}
	\begin{aligned}
	\mathcal{I}_1(p,q)=p^{n+1}+&\frac12 \delta t\Delta q^{n+1}+{\frac14} \delta t\big[((p^2+q^2)^{\mu}q)^{n+1}+[(p^2+q^2)^\mu]^{n+1}q^n\\
	&+[(p^2+q^2)^{\mu}]^nq^{n+1}\big]-p^n+\frac12 \delta t\Delta q^n+{\frac14} \delta t[(p^2+q^2)^{\mu}q]^n=0,\\
	\mathcal{I}_2(p,q)=q^{n+1}&-\frac12 \delta t\Delta p^{n+1}-{\frac14} \delta t\big[((p^2+q^2)^\mu p)^{n+1}+[(p^2+q^2)^\mu]^{n+1}p^n\\
	&+[(p^2+q^2)^\mu]^np^{n+1}\big]-q^n-\frac12 \delta t \Delta p^n-{\frac14} \delta t[(p^2+q^2)^\mu p]^n=0.
	\end{aligned}
	\end{align}
	One can parameterize $p$ and $q$ by a neural network and approximate the solution under the total residual loss $\mathcal{I}=\mathcal{I}_1+\mathcal{I}_2$.
	In this approach, we consider again those solutions that correspond to the Gaussian wave packets as initial conditions to build the datasets that will be  fed into FCNNs and ResNets in high dimensions.  
}
% {\red what is this for?}
%\end{remark}

\subsection{Optimization}
Formulated as an optimization problem, the parameters in the network can be obtained by using the stochastic or batch optimized algorithms, applied to the expected or empirical risks for \cref{eq: training} and \cref{eq: training-simo}. For a comparison of these methods, one can refer to \cite{bottou2018optimization}. The prototypical stochastic optimization method is the stochastic gradient descent method in \cite{robbins1951stochastic}, which, in the context of minimizing $\Lm(U_0,\Wm)$, with $\Wm_0$ initialized by \cite{he2015delving}, is defined by
\begin{align}
\Wm_{k+1} \leftarrow \Wm_{k}-\alpha_k \nabla_{\Wm_{i_k}} \Lm(U_0,\Wm_k),
\end{align}
for all $k\in \mathbb{N}$. The index $i_k$ is chosen randomly and $\alpha_k$ is a positive stepsize known as the learning rate. Each epoch of this method is thus very cheap, involving only the computation of the gradient {$\nabla_{\Wm_{i_k}}\Lm(U_0,\Wm_k)$} corresponding to one sample. In many cases, a batch approach is a more natural fit.
% more nature and well-known idea. The simplest such method in this class is the steepest descent algorithm (also called full-gradient method) which can be defined by
%\begin{align}
%    \Wm_{k+1} \leftarrow \Wm_k-\frac{\alpha_k}{N}\sum_{i=1}^N \nabla %\Lm_i(U_0,\Wm_k), 
%\end{align}
%such approach is more expensive than computing the step $-\alpha_k \nabla \Lm_{i_k}(U_0,\Wm_k)$ in SG, though one may expect a better step is computed when all samples are considered in an epoch. 
In this paper, we employ the \textit{Adam} method  \cite{kingma2014adam} at the beginning of the training process. {The convergence will be further improved by using the Broyden-Fletcher-Goldfarb-Shanno (L-BFGS) \cite{Nocedal1999NumericalO} method.}   

\section{Numerical Experiments}
\label{sec:experiments}

In this section, we present numerical examples to test the effectiveness of the neural network representation \eqref{eq: nnet0}. Extensive tests are performed to study the accuracy of the approximation and examine extrapolations by the neural networks. The training samples for the first two examples are based on  analytical solutions of the wave equation \eqref{wave-eq}, with the first example in  {1D} and the second example in  {3D}. We also extend the numerical test to wave equations in { 8 dimension} where the wave propagation occurs mainly in two dimensions. For the third example, we consider the linear Schr\"odinger equation \eqref{Schrodinger-1} and build the training set from analytical solutions. In the remaining { four} examples, we test our method for the cubic Schr\"odinger equations with solutions computed numerically {and the nonlinear Schr\"odinger equation with data generated from PINNs. }

\begin{example}[The 1D wave equation]
	
	Here we first consider the wave equation \cref{wave-eq} in 1D. The training sets are gathered by \cref{eq: dataset-single} and \cref{eq: dataset-multi}. 
	In the numerical experiments, we take the neural network with $D=5$ and $m_1= 2N_x$, $m_2=N_x$, $M=100$,
	% \begin{align*}
	% 	D=\{2N_x,100,100,100,100,N_x\}, 
	% \end{align*}
	both for the FCNNs and ResNets, in the latter case, we choose ResNets with two residual blocks, each block with $D=2$, $M=100$ and a skip connection so that the number of parameters of both networks are the same. Meanwhile, we take $N_x=201$, which is the number of grid points for both the training and testing samples in the spatial domain $[-8,8]$. The exact training samples are specified by sets $\mathcal{K}$ and $\Sigma$. We choose $\mathcal{K}=\{1,2,\cdots,10\}$. For $\Sigma$, we consider two types of selections: a set with linear spacing $L=\{0.8,0.9,1,1.1,1.2,1.3\}$, 
	and a set with exponential grid $E=\{h,2h,2^2h,2^3h,2^4h,2^5h\}$ where the spacing is doubled each time.
	We train the networks for $20000$ epochs { with $\lambda=1/16$}.

	After the parameters in the network are determined,  the performance of the network approximation is tested on solutions with  the following initial conditions,
	\begin{align}
	u^\i(x,0)&= \exp\left(-x^2\right)\cos(6x),\; u_t^\i(x,0)=\exp\left(-x^2\right)\left[2x\cos(6x)+6\sin(6x)\right],  \label{wave-initial-1}\\
	u^\ii(x,0)&= \exp\left(-x^2/1.5\right)\cos(6.5x),\; u_t^\ii(x,0)=\exp\left(-x^2/1.5\right)\left[\frac{x}{0.75}\cos(6.5x)+6.5\sin(6.5x)\right], \label{wave-initial-2}\\
	u^\iii(x,0)&=\sech(x)\cos(2x),\; u_t^\iii(x,0)=\sech(x)\left[\tanh(x)\cos(2x)+2\sin(2x)\right],\label{wave-initial-3}\\
	u^\iv(x,0)&= \exp\left(-x^2\right)\cos(\tilde{k}x),\; u_t^\iv(x,0)=\exp\left(-x^2\right)\left[2x\cos(\tilde{k}x)+\tilde{k}\sin(\tilde{k}x)\right],\label{wave-initial-extra-1}
	\end{align}
	with $\sech(x)=2/(\exp(x)+\exp(-x))$. 
	
	These initial conditions are selected based on the following rationale: We notice that  \cref{wave-initial-1} is of the same type of initial condition as those in the training sets presented in \cref{wave-initial-train-u} and \cref{wave-initial-train-ut}. It can be used to verify the training procedure.  The initial conditions in \cref{wave-initial-2} has a similar function form as those in the training set, but the wave number $k$ and the width $\sigma$ do not belong to  $\mathcal{K}$ and $\Sigma$. In view of the selection of $\mathcal{K}$ and $\Sigma$, this can be interpreted as an interpolation in terms of the wave number, but an {\it extrapolation} in terms of the width parameter.  The initial condition \cref{wave-initial-3} is outside of training sets in the sense that the function form is completely different.    For the last initial condition \cref{wave-initial-extra-1},  the wave number  $\tilde{k}$ will be selected as { $\tilde{k}=10.025 \;\textrm{and}~ 10.05$}
	to examine the extrapolation error.
	
	% to examine the extrapolation accuracy. For example, we can take the initial conditions as training sets with $\mathcal{K}=\{1,2,\cdots,10\}$, $\Sigma=\{0.8,0.9,1,1.1,1.2,1.3\}$ and predict the solutions with $\tilde{k}=10.025,10.05,10.1,10.2$ in a short range. The exact solution to \cref{wave-eq} can be found using the d'Alembert formula.
	%, while 

	%Note that the exact solution to \cref{wave-eq} given initial conditions
	%\begin{align*}
	%u(x,0)=\sech(x)\cos(kx),\quad u_t(x,0)=\sech(x)[tanh(x)\cos(kx)+k\sin(kx)]
	%\end{align*}
	%is drawn by
	%\begin{align*}
	%u(x,t)=\sech(x-t)\cos(k(x-t)).
	%\end{align*}

	Thanks to the availability of the exact solution, given by the d'Alembert formula, we can  quantify the error. Specifically, we define the relative error to be
	\begin{align}\label{eq-relative}
	\mathcal{R} =\frac{\|u^{\textrm{DNN}}-u^{\textrm{REF}}\|_2}{\|u^{\textrm{REF}}\|_2},
	\end{align}
	{ where notation $u^{\text{REF}}$ denotes the reference solution which can be taken by the exact solution, a numerically computed solution, or an approximated solution from PINNs.}
	
	\medskip
	
	The results, in terms of the relative error of the solutions at time $T=2$, are shown in \cref{tab-2-1}, where we collected the results for the solutions that correspond to the initial conditions \cref{wave-initial-1}, \cref{wave-initial-2} and \cref{wave-initial-3}, respectively.  The results for the first initial condition is hardly surprising, since the initial condition  \eqref{wave-initial-1} is very similar to those in the training set. But our numerical experiments suggest that this method also yields reasonable accuracy for the initial conditions \cref{wave-initial-2}  {and \cref{wave-initial-3}} that are not the type in the training set.
	We also observe that in most cases the ResNets yield slightly better accuracy. In this case, using the  {$\text{tanh}$ function with $L$ or $E|_{h=0.5}$} yields the best result. For the third case \cref{wave-initial-3}, we observe the choice of $\Sigma$ with linear spacing produces poor results, and it seems important to have a larger range of width parameters in the training set. The results also indicate   that 
	the choice of the activation function plays a role. { For example, with the choice of the FCNNs and ResNets, the $\text{sigmoid}$ function yields slightly worse results, while for the $\text{tanh}$ function, the accuracy is much better.  } 
	
	\begin{table}[ht]
		\centering
		\begin{tabular}{c|c|cc|cc|cc}
			\toprule   
			{}  & {Width}  & {$\mathcal{R}(\times 10^{-3})$}& &{$\mathcal{R}(\times 10^{-3})$} & &{$\mathcal{R}(\times 10^{-2})$} \\
			\cmidrule{2-2} 
			\cmidrule{3-4} 
			\cmidrule{5-6} 
			\cmidrule{7-8} 
			% 			\cmidrule{4-4} 
			$\phi$ & $\Sigma$ & $u_{\textrm{F}}$ & $u_{\textrm{R}}$   & $u_{\textrm{F}}$ & $u_{\textrm{R}}$ & $u_{\textrm{F}}$ & $u_{\textrm{R}}$  \\ 
			\midrule
			$\text{relu}(x)$ & $E|_{h=0.1}$ &{\ $3.03$}  & {\ $1.36$}  &\  $8.77$ &\  $2.92$ & \ $4.92$& \ $3.82$\\
			& $E|_{h=0.5}$ &{\ $3.15$} & {\ $2.50$}  & \ $7.84$&\ $8.13$ &\ $2.71$ &\ $2.21$ \\ 
			& $E|_{h=1}$ &{\ $3.91$} & {\ $3.60$} &\ $6.06$ &\ $14.4$&\ $3.09$ &\ $1.88$ \\ 
			& $L$ & {\ $2.31$} & {\ $0.985$} &\  $49.3$&\ $33.2$& \ $26.7$&\ $25.3$ \\ 
			\midrule
			$\tanh(x)$ & $E|_{h=0.1}$ & {\ $2.31$}  & {\ $1.14$} &\ $2.40$ &\ $2.05$& \ $4.26$&\ $3.43$  \\
			& $E|_{h=0.5}$ &{\ $2.44$}  & {\ $1.47$} & \ $2.10$& \ $1.43$&\ $0.651$ &\ $0.373$\\ 
			& $E|_{h=1}$ & {\ $3.61$} & {\ $2.21$} &\  $2.92$&\ $1.81$ &\ $0.690$ & \ $0.347$\\ 
			& $L$ & {\ $0.797$} & {\ $0.550$} & \ $4.61$&\ $3.53$& \ $12.3$& \ $11.8$\\ 
			\midrule
			$\text{sigmoid}(x)$ & $E|_{h=0.1}$ &  {\ $10.4$} & {\ $2.95$} &\ $2.77$ &\ $4.12$ &\ $7.39$ &\ $3.21$ \\
			& $E|_{h=0.5}$ & {\ $9.82$} & {\ $3.07$} &\  $14.9$&\ $3.55$&\ $5.95$ &\ $0.642$ \\ 
			& $E|_{h=1}$ & {\ $1.32$} & {\ $5.97$} & \ $26.0$&\ $5.69$ & \ $6.68$&\ $0.802$ \\ 
			& $L$ &{\ $6.41$}  &{\ $1.69$} & \ $29.6$&\ $12.5$ & \ $21.9$& \ $13.2$\\ 
			\midrule
			$\text{elu}(x)$ & $E|_{h=0.1}$ &  {\ $1.60$} & {\ $1.20$} &\ $3.61$ &\ $2.47$&\ $3.82$ &\ $3.50$  \\
			& $E|_{h=0.5}$ &{\ $2.16$} & {\ $1.54$} &\ $2.97$ &\ $2.03$&\ $1.00$ &\ $0.443$  \\ 
			& $E|_{h=1}$ &  {\ $2.79$} & {\ $2.01$} &\ $2.88$ &\ $1.91$& \ $0.811$& \ $0.456$ \\ 
			& $L$ &  {\ $1.05$} & {\ $0.798$} &\ $6.24$ &\ $4.26$ &\ $12.6$ &\ $12.2$ \\ 
			
			\bottomrule
		\end{tabular}

		\caption{The approximation error for various choices of activation functions and band width of wave packets for the 1D wave equation \eqref{wave-eq} with FCNNs and ResNets for the initial conditions \cref{wave-initial-1} (Left), \cref{wave-initial-2} (Middle) and \cref{wave-initial-3} (Right). In the training sets, the wave numbers are chosen from  $\mathcal{K}=\{1,2,\cdots,10\}$ { and the rescaling parameter $\lambda=1/16$.}}\label{tab-2-1}
	\end{table}

	The relative error of solutions at time $T=2$ from the initial condition \cref{wave-initial-extra-1}
	with various choices of $\tilde{k}$ outside  the training set, is shown in \cref{tab-extra-1}. The training sets are constructed with $\mathcal{K}=\{1,2,\cdots,10\}$ and $\Sigma=L$. For $\tilde{k}$ from $10.025$ to $10.05$ the accuracy  is reasonable. But we do observe that it deteriorates as   $\tilde{k}$ moves further away from $\mathcal{K}$. Interestingly, for the ResNet with activation functions {$\text{tanh}(x)$} and $\text{elu}(x)$, the error grows much more slowly. 
	
	\begin{table}[h]
		\centering
		\begin{tabular}{lcccccc}
			\toprule   
			{}  & $\mathcal{K}$  & $\mathcal{R}(\times 10^{-3})$ & &   $\mathcal{K}$  & $\mathcal{R}(\times 10^{-3})$& \\
			\cmidrule{2-2} 
			\cmidrule{3-4} 
			\cmidrule{5-5} 
			\cmidrule{6-7} 
			% 			\cmidrule{4-4}
			$\phi$ & $\tilde{k}$ & $u_{\textrm{F}}$ & $u_{\textrm{R}}$ & $\tilde{k}$ & $u_{\textrm{F}}$ & $u_{\textrm{R}}$ \\ 
			\midrule
			$\text{relu}(x)$ & $10.025$ & {\ $7.26$}& {\ $5.82$}   
			& $10.05$ & {\ $14.1$} & {\ $11.9$}  \\ 
			%		& $10.1$ & {\ $28.4$} & {\ $25.1$} \\ 
			%		& $10.2$ &{\ $58.9$} & {\ $55.3$} \\ 
			\midrule
			$\tanh(x)$ & $10.025$ & {\ $1.76$} & {\ $1.54$}  
			& $10.05$ & {\ $1.87$} & {\ $1.64$} \\ 
			%		& $10.1$ & {\ $2.16$} & {\ $1.86$}  \\ 
			%		& $10.2$ & {\ $2.95$} & {\ $2.45$}  \\ 
			\midrule
			$\text{sigmoid}(x)$ & $10.025$ & {\ $6.71$} & {\ $4.13$}  
			& $10.05$ & {\ $12.5$} & {\ $4.36$}  \\ 
			%		& $10.1$ & {\ $25.8$} & {\ $4.86$} \\ 
			%		& $10.2$ & {\ $55.2$} & {\ $6.01$}  \\ 
			\midrule
			$\text{elu}(x)$ & $10.025$ & {\ $1.78$} & {\ $1.73$}  
			& $10.05$ & {\ $2.11$} & {\ $1.88$} \\ 
			%		& $10.1$ & {\ $3.17$} & {\ $2.36$}  \\ 
			%		& $10.2$ & {\ $6.13$} & {\ $3.85$}  \\ 
			
			\bottomrule
		\end{tabular}
		
		\caption{Extrapolating the wave number. The approximation of   the 1D wave equation with the initial condition \cref{wave-initial-extra-1} using different activation functions and band width parameters  with FCNNs and ResNets. Choosing $\mathcal{K}=\{1,2,\cdots,10\}$, $\Sigma=L$ and $\lambda=1/16$.}\label{tab-extra-1}
	\end{table}
\end{example}

\begin{example}[High dimensional wave equations]
	In high dimensions, in general, the wave modes are represented by many wave numbers. Here we consider a special case where the variation of the solutions of   \cref{wave-eq} is mainly in the first two dimensions. To this end, we choose training samples specified by $\mathcal{K}=(k_1,k_2,k_3,\cdots,k_d)$ with $k_1$, $k_2$ both in $\{1,2,\cdots,5\}$ and $k_i=1,\;i=3,\cdots,d$, which indicates that the wave propagation is mostly restricted to they $xy$-plane. We also choose $\text{tanh}(x)$ and { $\Sigma=L$}.
	
	Starting with the initial condition  \cref{wave-initial-train-u}  {and \cref{wave-initial-train-ut}}, we use neural network to represent \cref{eq: nnet0} with $\k=(2,2,1,\cdots,1)$ and $2\sigma^2=1$. We take the neural network with $D=5$, $m_1=2N_x N_y$, $m_2=N_x N_y$ and $M=100$ of FCNNs.  We also take $N_x=N_y=64$ as the number of grid points for both training and testing in domain $\Omega=[-4,4]\times [-4,4]$. We train the networks of { $40000$} epochs. We increase the dimension from $3$D to $ 8$D. The solution at the cross section with the $xy-$plane  at $T=0.5$ are shown in  \cref{fig-wave-highD}. The relative error are { $4.48\times 10^{-3}$ (3D) and $1.02\times 10^{-2}$ (8D)}. One can observe that the error grows as  the dimension increases. 
	
	\begin{figure}[htbp]
		\centering
		%	\subfloat{\textbf{(3D)}}
		\subfloat{\includegraphics[width=1.6in]{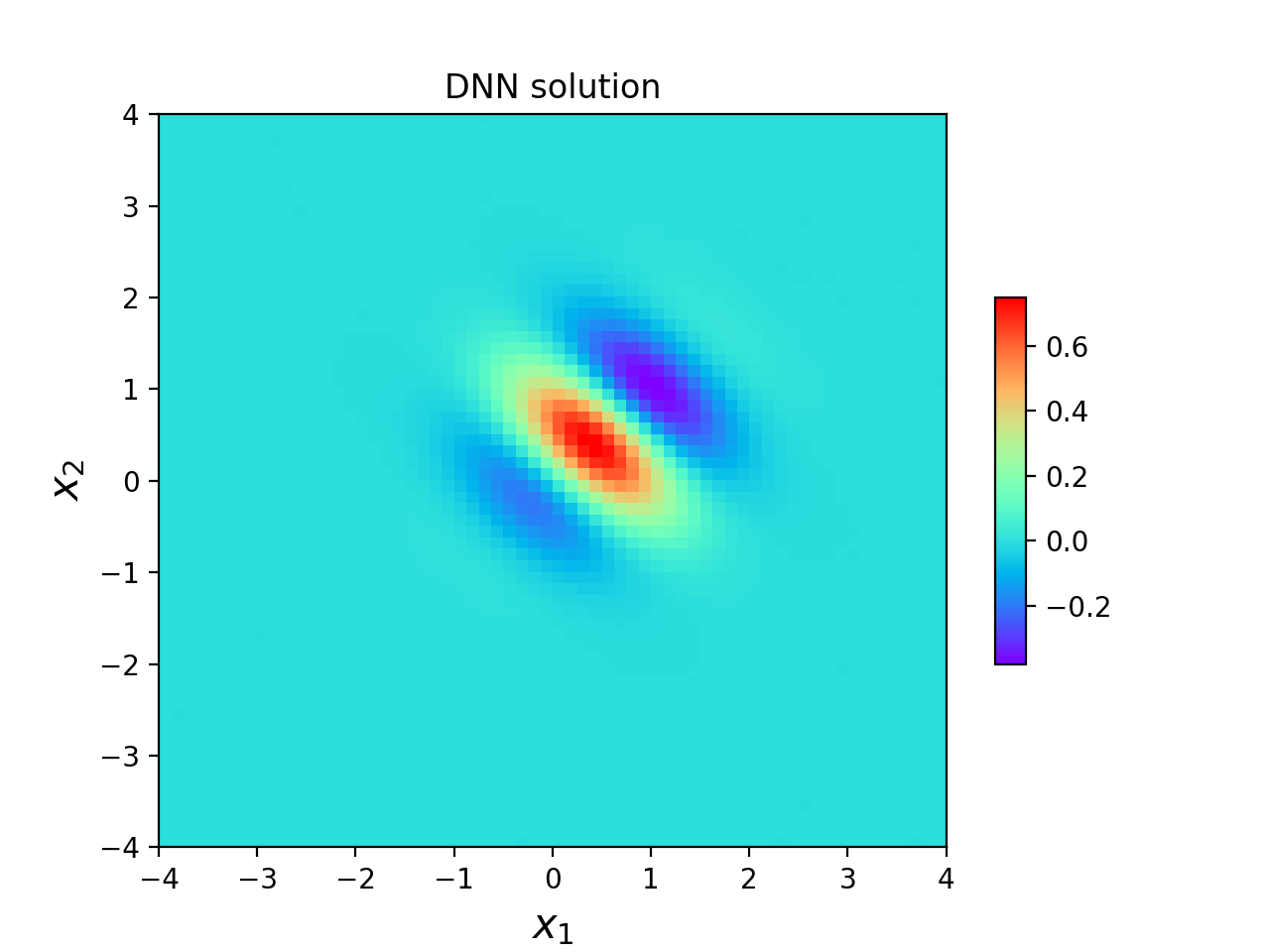}}
		\subfloat{\includegraphics[width=1.6in]{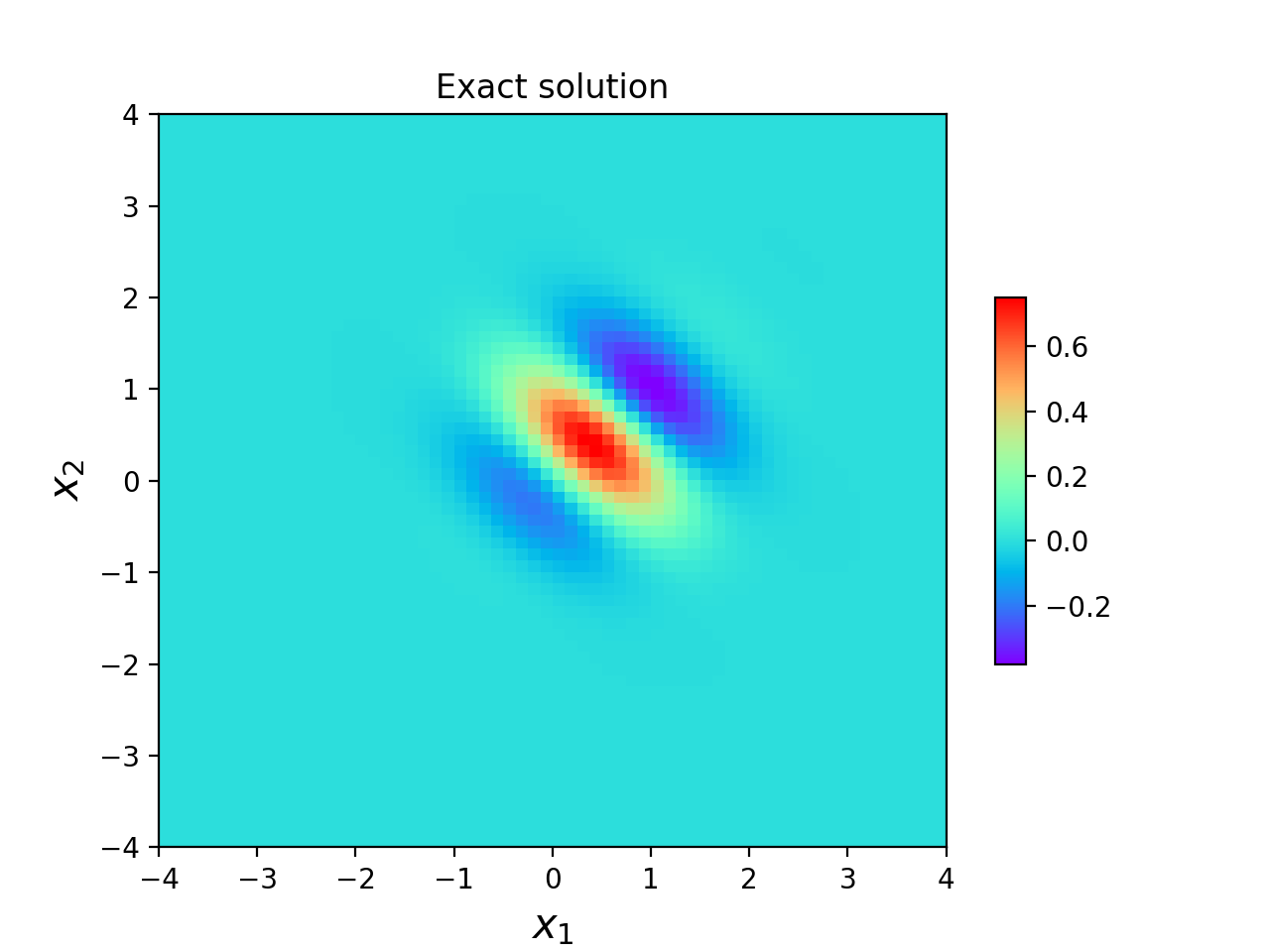}}
		%	\subfloat{\includegraphics[width=2in]{wave_err_3D_v1_V1.png}}
		%	\hspace{0.1in}
		%	\subfloat{\textbf{(8D)}}
		\subfloat{\includegraphics[width=1.6in]{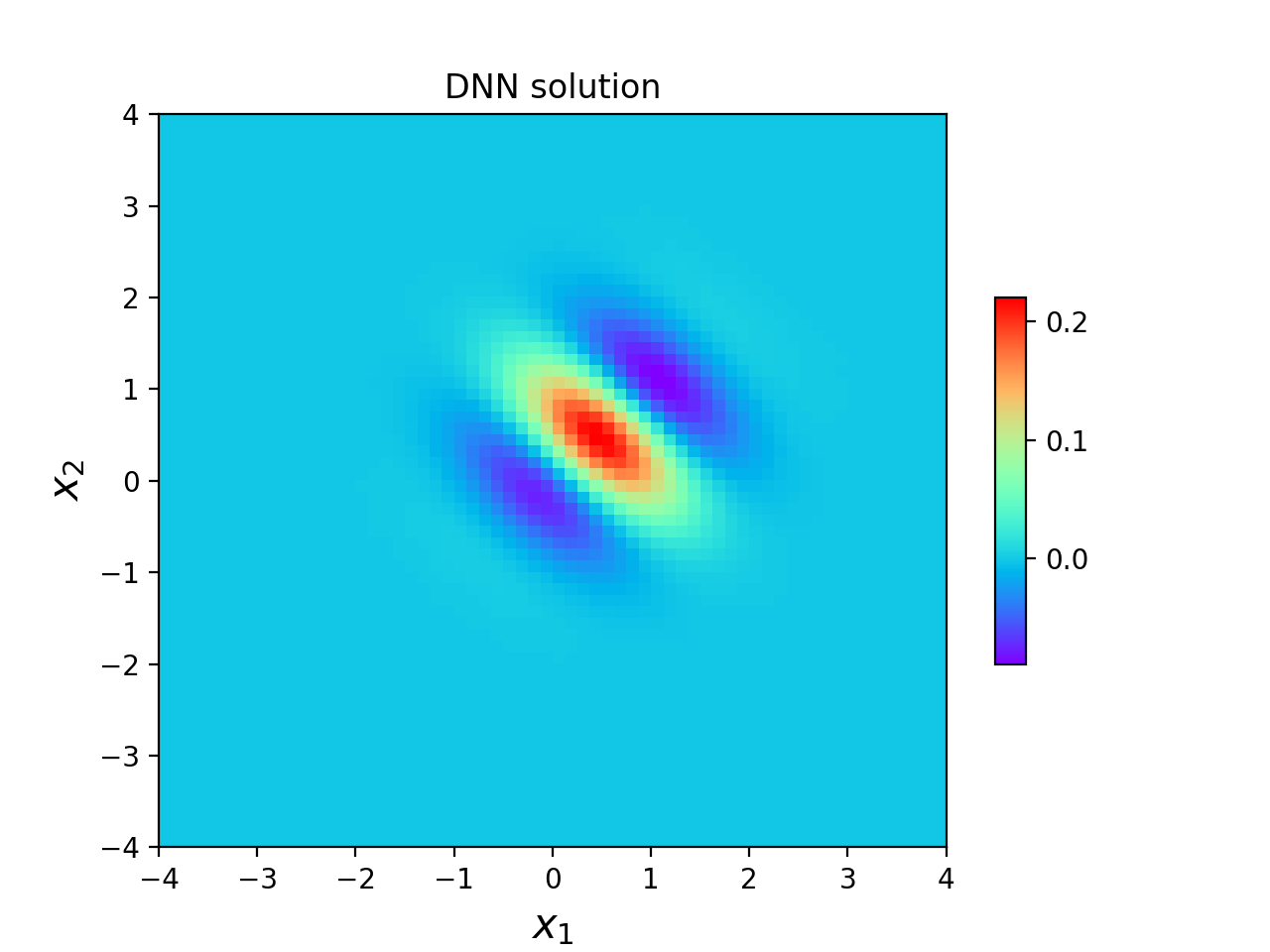}}
		\subfloat{\includegraphics[width=1.6in]{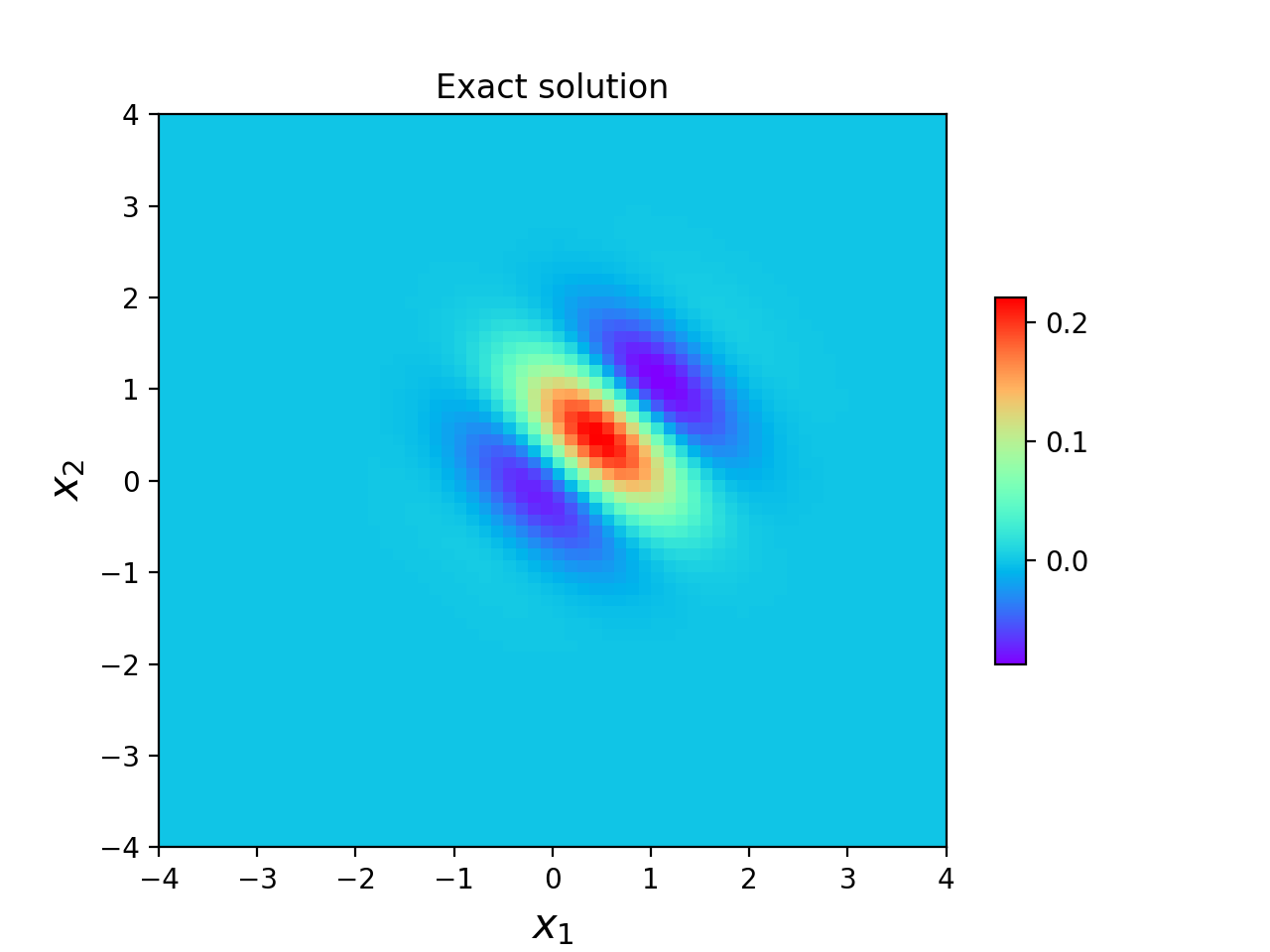}}
		%	\subfloat{\includegraphics[width=2in]{wave_err_8D_v1_V1.png}}
		%	\hspace{0.1in}
		%	\subfloat{\textbf{(16D)}}
		%	\subfloat{\includegraphics[width=2in]{wave_dnn_sol_16D_v1.png}}
		%	\subfloat{\includegraphics[width=2in]{wave_exact_sol_16D_v1.png}}
		%	\subfloat{\includegraphics[width=2in]{wave_err_16D_v1.png}}
		\caption{Solution of high-dimensional wave equation by DNN (Left) {\  with $\text{tanh}(x)$ up to $40000$ iterations}, exact solution (Right) evaluated at time $T=0.5$ to high dimensional wave equation. Left two panel: $d=3$; Right two panel: $d=8$.}\label{fig-wave-highD}
	\end{figure}

\end{example}

\begin{example}[The 1D linear Schr\"odinger equation]
	We consider the linear Schr\"odinger equation \cref{S-2} in the 1D case. The training sets are constructed from \cref{eq-p1} and \cref{eq-1D}. We consider the following initial conditions, 
	% {\color{red}what the value of $\tilde{k}$?}
	\begin{align}
	u_0^{\i}(x)&=\exp\left[-(\frac{x}{\sqrt{2}})^2 + i 3\sqrt{2}x\right],\label{initial-Schr-1}\\
	u_0^{\ii}(x)&=\exp\left[-(\frac{x}{\sqrt{2}})^2/1.2 + i 3.25\sqrt{2}x\right], \label{initial-2-1}\\
	u_0^{\iii }(x)&=\exp\left[-(\frac{x}{\sqrt{2}})^2 + i\tilde{k}\frac{x}{\sqrt{2}}\right]. \label{initial-Sch-extra-k}
	%\psi_0^{\textrm{re}}(x)&= sech(\frac{x}{\sqrt{2}})\cos(6\frac{x}{\sqrt{2}}), \quad 	\psi_0^{\textrm{im}}(x) = sech(\frac{x}{\sqrt{2}})\sin(6\frac{x}{\sqrt{2}}).\label{initial-Schr-2}
	\end{align}
	In the third case,  $\tilde{k}$ is to be selected to examine the extrapolation error. %For example, the solutions can be predicted with $\tilde{k}=6.025$, or $6.05, 6.1, 6.2$ in terms of the initial conditions as training sets with $\mathcal{K}=\{1,2,\cdots,6\}$, and $\Sigma=\{0.8,0.9,1,1.1,1.2,1.3\}$.
	
	% \begin{align}\label{initial-2-1}
	% \psi_0^{\textrm{re}}(x)=\exp(-x^2/1.2)\cos(6.5x), \quad 	\psi_0^{\textrm{im}}(x) = \exp(-x^2/1.2)\sin(6.5x).
	% \end{align}
	We consider the neural network with $D=5$, $m_1=N_x$, $m_2=N_x$ and $M=100$ of FCNNs and two residual blocks for ResNets, each block with $D=2$, $M=100$.  $N_x=201$ is the  number of grid points for both training and testing in the domain $\Omega=[-8,8]$. The training samples are specified by $\mathcal{K}$ and $\Sigma$. We train the networks for $20000$ epochs. 
	
	We first consider initial conditions \cref{initial-Schr-1} with $\lambda=1/4$ and \cref{initial-2-1} with $\lambda=1/2$ as input, and the corresponding solutions  at a single time  instance $T=0.2$ as the output. The results are shown in \cref{tab-2}. For the initial condition \cref{initial-Schr-1} the best result is obtained by using the 
	$\text{relu}$ function for both FCNN and ResNet. For the initial condition \cref{initial-2-1},  the accuracy is not as satisfactory as the previous test, especially when the 
	$\text{relu}$ function is used. Another observation is that the result is quite sensitive to the selection of $\Sigma$ for the training set.

	\begin{table}[htbp]
		\centering
		\begin{tabular}{@{\extracolsep{4pt}}c|cc|cc|ccccc}
			\toprule   
			{Width}  & \multicolumn{2}{c}{$\mathcal{R}({\tiny\times 10^{-3}}), \; \lambda_1$}& \multicolumn{2}{c}{$\mathcal{R}({\tiny\times 10^{-2}}), \; \lambda_2$}&	{Width}  & \multicolumn{2}{c}{$\mathcal{R}({\tiny\times 10^{-3}}), \; \lambda_1$}& \multicolumn{2}{c}{$\mathcal{R}({\tiny\times 10^{-2}}), \; \lambda_2$}\\
			\cmidrule{1-1} 
			\cmidrule{2-3} 
			\cmidrule{4-5} 
			\cmidrule{6-6} 
			\cmidrule{7-8} 
			\cmidrule{9-10} 
			%		\cmidrule{5-7} 
			%  \cmidrule{4-6} 
			$\Sigma$ & $\rho_{\textrm{F}}$ & $\rho_{\textrm{R}}$& $\rho_{\textrm{F}}$ & $\rho_{\textrm{R}}$&$\Sigma$ & $\rho_{\textrm{F}}$ & $\rho_{\textrm{R}}$& $\rho_{\textrm{F}}$ & $\rho_{\textrm{R}}$\\ 
			\midrule
			$E|_{h=0.1}$ & {$9.19$} & {$10.7$}& $4.72$&$3.86$&$E|_{h=1}$ &  {$2.24$} & $2.24$ &$1.20$ &$4.28$ \\
			$E|_{h=0.5}$  & {$2.38$} & $2.14$ &$2.07$ &$5.62$ &$L$ & {$2.96$}  & $1.20$ &$2.90$ &$1.81$\\ 
			\midrule
			$E|_{h=0.1}$& {$4.88$}&$9.53$& $1.17$&
			$0.94$&	$E|_{h=1}$  &{$3.02$} & $5.27$&$0.82$ &$0.54$ \\
			$E|_{h=0.5}$ &{$7.60$} &$8.63$  &$0.63$ &$0.59$ &$L$ &  {$3.47$} & $3.89$ &$1.93$ &$0.75$\\ 
			\midrule
			$E|_{h=0.1}$ & {$10.8$} &  $20.6$& $7.55$&$1.98$ &$E|_{h=1}$ &  {$4.91$} &$7.00$&$1.08$ &$0.89$ \\
			$E|_{h=0.5}$ &{$11.2$} & $20.9$ &$1.65$ & $1.27$&$L$ & {$3.34$}  &$5.09$ &$2.15$ &$1.13$ \\ 
			\midrule
			$E|_{h=0.1}$  & {$7.02$}  &$9.47$ &$1.62$ &$1.00$&$E|_{h=1}$ & {$2.88$} &$3.05$ & $0.84$& $0.74$ \\
			$E|_{h=0.5}$ & {$6.22$} &$6.83$&$0.58$ &$0.69$&$L$ & {$3.36$} &$2.75$& $1.28$&$1.01$ \\ 
			
			\bottomrule
		\end{tabular}
		%	\mskip
		\caption{Approximation error {up to $T=0.2$} for the 1D linear Schr\"odinger equation with  initial conditions \cref{initial-Schr-1} and \cref{initial-2-1}  using FCNNs and ResNets.  From top to bottom: Different choices of  activation functions $\textrm{relu}(x)$, $\tanh(x)$, $\textrm{sigmoid}(x)$ and $\textrm{elu}(x)$, and  various choices of $\Sigma$. The wave numbers are drawn from $\mathcal{K}=\{1,2,\cdots,10\}$ to build the training set. { Taking the scaling parameter $\lambda=\lambda_i,\;i=1,2$ with $\lambda_1=1/4$, $\lambda_2=1/2$.}}\label{tab-2}
	\end{table}
	
	Next we discuss the results from the extrapolation. In this context, an extrapolation can be interpreted in terms of the wave number  $\tilde{k}$ 
	in the initial condition \cref{initial-Sch-extra-k}, or in terms of predicting solutions at time instances that are beyond the training period. In the former case, 
	we consider a training set determined by $\mathcal{K}=\{1,2,\cdots,6\}$, $\Sigma=L$, and {$\lambda=1/2$} and then we pick an initial condition  \cref{initial-Sch-extra-k}, where { $\tilde{k} \in \{6.025,6.05\}$} The results 
	are summarized in \cref{tab-extra-K}. {The best result comes from the FCNN with $\text{elu}(x)$.} One can see that the error is reasonable for wave numbers in this range, but  in all cases the error increases as $\tilde{k}$  
	moves further away from $\mathcal{K}$. {To visualize the extrapolation in wave number, we take the FCNN with $\elu(x)$ for the best performance which is shown in \cref{fig-extra-K-S}.} In the latter case, we consider the solutions 
	with the initial condition \cref{initial-Schr-1}.  The training set  consists of solutions  $(x,t)\in [-8,8]\times [0,0.6]$ with $N_x=201$, {$N_t=11$}. Then we test the solution  at time $T$ from $0.625$ to {$0.65$},  and the results are presented in \cref{tab-extra-schr-t}. The best result comes from the FCNN using the {$\elu(x)$} activation function. But the error grows with $T$, indicating that the accuracy of the extrapolation can only be guaranteed for a finite time period. To visualize the evolution in time, we take the FCNN with the best performance, together  with the activation function {$\elu(x)$}. The results are  shown in \cref{fig-extra-T-1}.  A good performance can be guaranteed in short time.  But we can see noticeable error  for longer times.
	
	% ---------- extrapolation in k --------
	\begin{table}[htbp]
		\centering
		\begin{tabular}{@{\extracolsep{4pt}}ccccccc}
			\toprule   
			$\mathcal{K}$  & \multicolumn{3}{c}{$\mathcal{R}(\times 10^{-3})$}\\
			\cmidrule{1-1} 
			\cmidrule{2-4} 
			\cmidrule{5-7} 
			%  \cmidrule{4-6} 
			$\tilde{k}$ & $p_{\textrm{F}}$ & $q_{\textrm{F}}$ & $\rho_{\textrm{F}}$ & $p_{\textrm{R}}$ & $q_{\textrm{R}}$ & $\rho_{\textrm{R}}$\\ 
			\midrule
			$6.025$ & $7.58$&$8.33$ &$5.43$  &$6.03$ &$6.32$ &$4.04$ \\
			$6.05$ & $14.3$& $16.2$ & $11.5$ &$11.4$ &$12.2$ &$7.70$ \\ 
			%		$6.1$ &\red$2.83$ & \red$3.28$ &\red $2.38$ &\red$2.33$ &\red$2.44$ &\red$1.62$ \\ 
			%		$6.2$ & \red$6.02$& \red$7.29$ & \red$5.43$ &\red $5.44$ &\red$5.00$ &\red$3.85$ \\ 
			\midrule
			$6.025$ &$4.78$ &$4.71$ &$3.62$ & $3.93$&$3.64$ &$3.22$ \\
			$6.05$ & $7.79$& $7.64$ &$4.80$  &$6.47$ &$556$ &$4.34$ \\ 
			%		$6.1$ & \red$1.49$& \red$1.46$& \red$0.798$&\red$1.26$ & \red$1.04$&\red$0.714$ \\ 
			%		$6.2$ & \red$3.17$ &\red$3.09$ &\red$1.58$ &\red$2.74$ &\red$2.24$ & \red$1.39$\\ 
			\midrule
			$6.025$ & $7.65$&$10.1$ & $6.28$ &$9.17$ & $6.76$&$9.83$  \\
			$6.05$ &$13.1$ &$18.8$ &$12.0$ &$11.6$ &$8.82$ &$10.7$ \\ 
			%		$6.1$ & \red$2.60$& \red$3.77$& \red$2.51$ &\red$1.73$ &\red$1.47$ & \red$1.29$\\ 
			%		$6.2$ & \red$5.49$&\red$7.88$ &\red$5.42$ &\red$3.15$ &\red$3.02$ &\red$1.87$ \\ 
			\midrule
			$6.025$ & $3.42$& $4.03$& $3.49$ & $3.78$&$3.92$ &$3.45$ \\
			$6.05$ &$5.20$ &$6.42$ &$4.71$  & $6.19$&$6.13$ &$4.81$ \\ 
			%		$6.1$ &\red$0.976$ &\red$1.22$ & \red$0.785$&\red$1.19$ &\red$1.15$  &\red$0.802$ \\ 
			%		$6.2$ &\red$2.10$ &\red$2.57$ &\red$1.56$  &\red$2.54$ &\red$2.47$ &\red$1.59$ \\ 
			
			\bottomrule
		\end{tabular}
		
		\caption{Extrapolation in the wave number at $T=0.2$ with FCNNs and ResNets for the initial condition \cref{initial-Schr-1}. From top to bottom: Different choices of the activation functions, $\textrm{relu}(x)$, $\tanh(x)$, $\textrm{sigmoid}(x)$ and $\textrm{elu}(x)$. Choosing $\mathcal{K}=\{1,2,\cdots,6\}$, $\Sigma=L$ {and $\lambda=1/2$}.}\label{tab-extra-K}
	\end{table}
	
	\begin{figure}[htbp]
		\centering
		\subfloat{\includegraphics[width=3.3in]{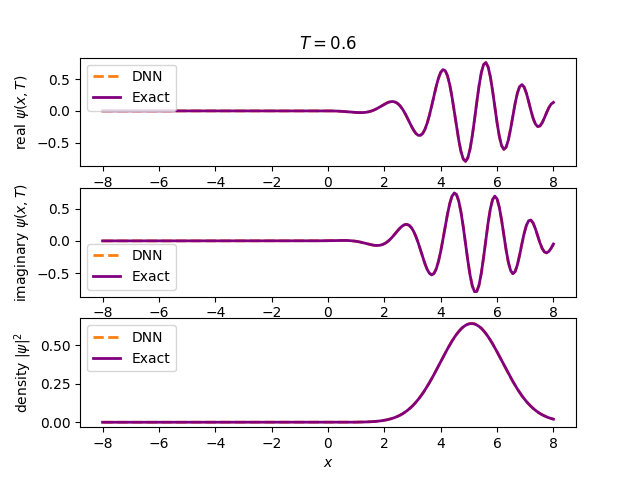}}
		\subfloat{\includegraphics[width=3.3in]{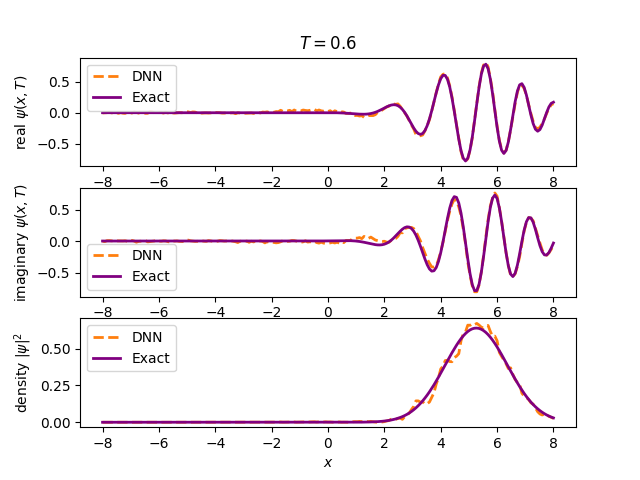}}
		%	\hspace{0.1in}
		%	\subfloat{\includegraphics[width=3in]{Schr_1D_time_extrapolation_0p6_T_0p675_v1_v1.png}}
		%	\subfloat{\includegraphics[width=3in]{Schr_1D_time_extrapolation_0p6_T_0p7_v1_v1.png}}
		\caption{Profile of the learned solution at $T=0.6$ for the 1D linear Schr\"odinger equation associated with the extrapolation in the wave number by using the activation $\text{elu}(x)$ of FCNNs. Choosing $\mathcal{K}=\{1,2,\cdots,6\}$, $\Sigma=L$ {and $\lambda=1/2$. Left: $k=6$; Right: $k=6.2$.}.}\label{fig-extra-K-S}
	\end{figure}
	
	% ------ extrapolation in time -------
	\begin{table}[htbp]
		\centering
		\begin{tabular}{@{\extracolsep{4pt}}ccccccc}
			\toprule   
			{Later time}  & \multicolumn{3}{c}{$\mathcal{R}(\times 10^{-2})$}\\
			\cmidrule{1-1} 
			\cmidrule{2-4} 
			\cmidrule{5-7} 
			%  \cmidrule{4-6} 
			$T$ & $p_{\textrm{F}}$ & $q_{\textrm{F}}$ & $\rho_{\textrm{F}}$ & $p_{\textrm{R}}$ & $q_{\textrm{R}}$ & $\rho_{\textrm{R}}$\\ 
			\midrule
			$0.625$ &$4.69$ & $5.01$ & $4.45$ &$9.45$ & $9.97$& $8.44$\\
			$0.65$ & $7.92$ & $9.68$ & $7.87$ &$18.7$ &$18.4$ & $14.3$\\ 
			%		$0.675$ &\blue$11.3$ & \blue$15.4$ &\blue$12.9$  &$73.4$ &$68.5$ & $7.14$\\ 
			%		$0.7$ & \red$16.0$ & \red$20.0$  &\red$14.1$  &$78.2$ &$73.5$ &$8.05$ \\ 
			\midrule
			$0.625$ &$1.46$ & $1.40$ &$1.31$ & $2.04$ & $2.33$&$1.56$ \\
			$0.65$ &$2.95$ & $3.08$ & $2.97$  &
			$3.36$ &$4.83$ & $2.86$\\ 
			%		$0.675$ &\red$5.84$ & \red$4.67$ & \red$6.23$ &$19.8$ &$28.0$ & $12.8$\\ 
			%		$0.7$ & \red$8.81$ & \red$6.69$ & \red$9.59$ &$24.7$ &$38.0$ & $13.9$ \\ 
			\midrule
			$0.625$ & $1.26$&$1.74$ &$1.65$  & $2.42$&$2.68$ &$2.45$  \\
			$0.65$ & $2.16$ &$2.63$ & $2.85$& $4.38$&$5.20$ &$4.43$ \\ 
			%		$0.675$ & $55.8$ & $44.5$ & $2.41$ &$27.1$ &$33.0$ & $7.67$\\ 
			%		$0.7$ & $59.7$ & $48.6$ &$2.57$  &$31.6$ &$38.6$ &$7.83$ \\ 
			\midrule
			$0.625$ &$1.31$ & $1.10$ & $1.36$  & $5.12$ & $4.07$ &$3.85$ \\
			$0.65$ & $2.11$ & $1.94$ &$2.43$  &$11.7$  &$8.19$ &$9.20$ \\ 
			%		$0.675$ & $15.0$& $12.7$&$8.20$  &$65.3$ &$62.6$  &$11.2$ \\ 
			%		$0.7$ & $17.5$& $13.7$ & $8.90$ & $68.5$ &$68.3$ &$14.0$ \\ 
			
			\bottomrule
		\end{tabular}
		
		\caption{Illustration of extrapolation in time with training sets in $[0,0.6]$. Different activation functions $\textrm{relu}(x)$, $\tanh(x)$, $\textrm{sigmoid}(x)$ and $\textrm{elu}(x)$ from the top row to the  bottom row for the 1D linear Schr\"odinger equation with the initial condition \cref{initial-Schr-1} as an input. Choosing $\mathcal{K}=\{1,2,\cdots,6\}$, $\Sigma=L$ {and $\lambda=1/2$}.}\label{tab-extra-schr-t}
	\end{table}
	
	\begin{figure}[htbp]
		\centering
		\subfloat{\includegraphics[width=3.3in]{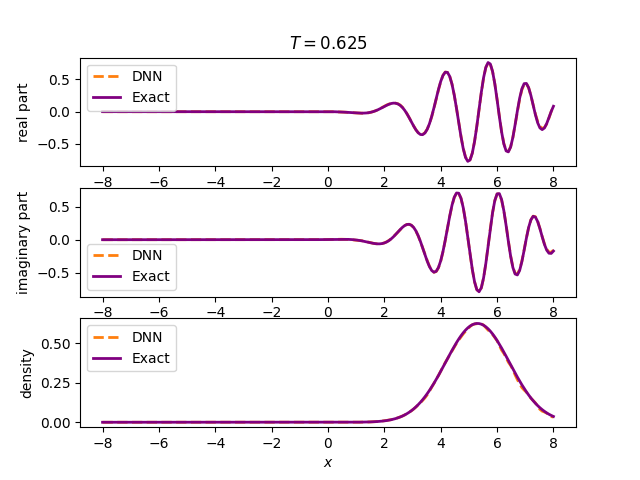}}
		\subfloat{\includegraphics[width=3.3in]{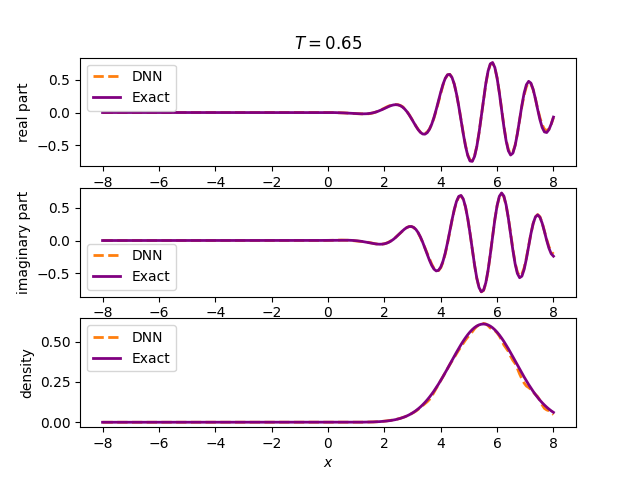}}
		%	\hspace{0.1in}
		%	\subfloat{\includegraphics[width=3in]{Schr_1D_time_extrapolation_0p6_T_0p675_v1_v1.png}}
		%	\subfloat{\includegraphics[width=3in]{Schr_1D_time_extrapolation_0p6_T_0p7_v1_v1.png}}
		\caption{Time evolution of the learned solution for the 1D linear Schr\"odinger equation associated with time extrapolation by using the activation $\text{elu}(x)$ of FCNNs. Choosing $\mathcal{K}=\{1,2,\cdots,6\}$, $\Sigma=L$ {and $\lambda=1/2$. Left: $T=0.625$; Right: $T=0.65$.}}\label{fig-extra-T-1}
	\end{figure}

	Next we test the accuracy of  a network trained using a dataset that consists of multiple snapshots of the solutions in the time interval $t\in[0,3]$ using uniform step size with $N_t=51$. The network is a FCNN with $D=5$, $m_1=N_x$, $m_2=N_t N_x$ and $M=100$, {activation function $\tanh(x)$}, $\mathcal{K}=\{1,2,\cdots,10\}$ and $\Sigma = \{0.8,0.9,1,1.1,1.2,1.3\}$.  Starting from the initial condition \cref{initial-Schr-1},  \cref{fig-4} shows the prediction by the FCNN, compared to the exact solution. The relative error for the density, the real part, the imaginary part of the wave function are given by {$8.900\times 10^{-3}$, $8.025\times 10^{-3}$, $7.606\times 10^{-3}$}, respectively.   Such examples appear frequently in testing an absorbing boundary condition \cites{arnold2003,jiang2004fast,wu2020absorbing}, and the main emphasis is usually on the reflection at the boundary. The results in   \cref{fig-4} suggest that the approximation by a FCNN exhibits an absorbing property that is similar to an absorbing boundary condition. 

	\begin{figure}[htbp]
		\centering
		\subfloat{\includegraphics[width=2in]{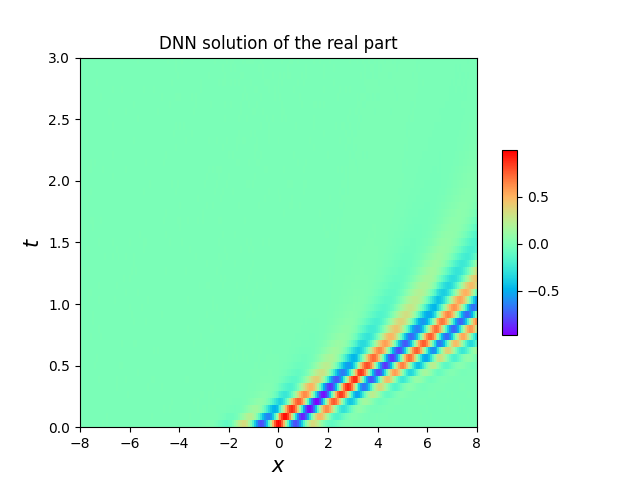}}
		\subfloat{\includegraphics[width=2in]{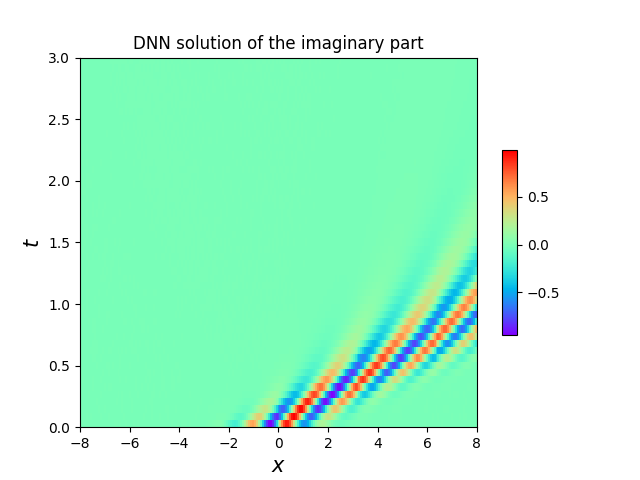}}
		\subfloat{\includegraphics[width=2in]{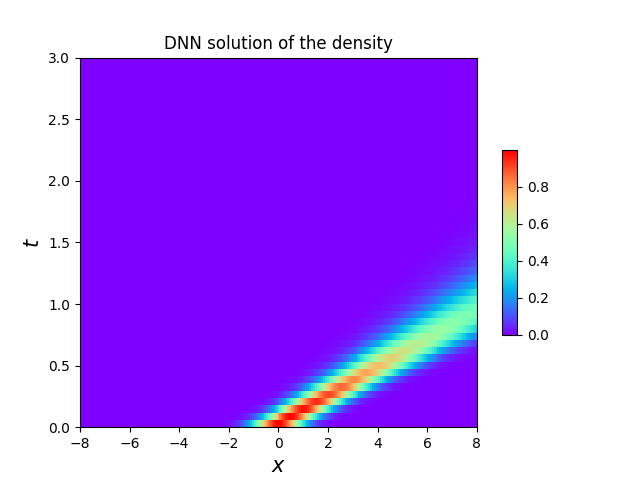}}
		\hspace{0.1in}
		\subfloat{\includegraphics[width=2in]{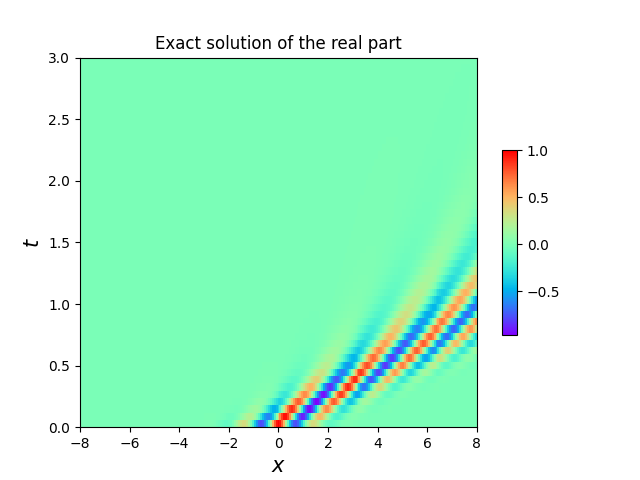}}
		%	\subfloat{\includegraphics[width=2in]{Linear_Schr_err_re_v2_t_V1.png}}
		\subfloat{\includegraphics[width=2in]{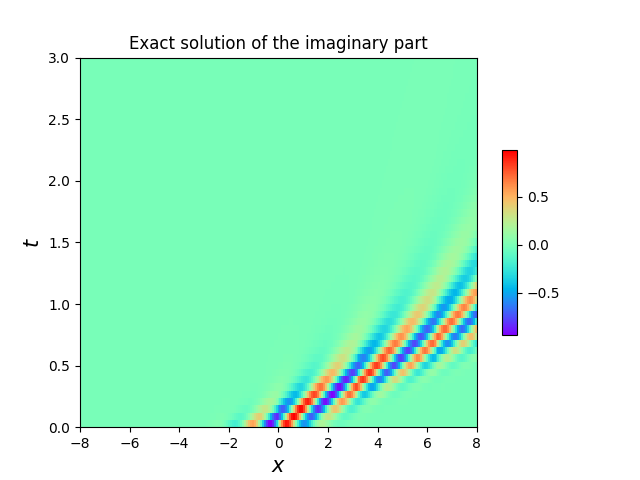}}
		%	\subfloat{\includegraphics[width=2in]{Linear_Schr_err_imag_v2_t_V1.png}}
		\subfloat{\includegraphics[width=2in]{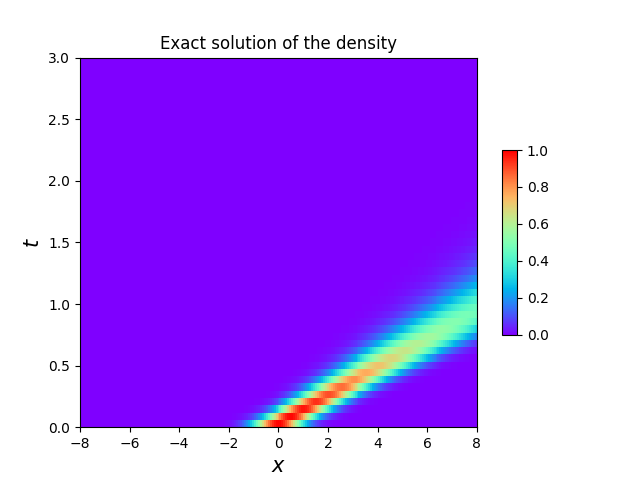}}
		%	\subfloat{\includegraphics[width=2in]{Linear_Schr_rho_err_v2_t_V1.png}}
		\caption{A wave packet propagating outside the domain { for 1D Schr\"odinger equation with \cref{initial-Schr-1}  as an input}. Top:  Prediction by a FCNN with $\text{tanh}(x)$; Bottom:  the exact solution. 
			% evaluated in time domain $t\in[0,3]$ for the initial condition \cref{initial-Schr-1}. 
			Left: the real part of wave packet solution; Middle: the imaginary part;  Right:  the electron density $\rho$.}\label{fig-4}
	\end{figure}
	
	%{In this framework, to see the hyper-parameter effect of choosing $\mathcal{K}$ and $\Sigma$, given the initial condition as
	%\begin{align}
	%    \psi_0^{\textrm{re}}(x)=\exp(-x^2)\cos(6x),\quad \psi_0^{\textrm{im}}(x)=\exp(-x^2)\sin(6x),
	%\end{align}
	%which is in the training set.
	%The setup for our neural network is taken number of layers $D=\{201,100,100,100,100,201\}$ and the $N=201$ number of grid points for both training and testing in domain $[-8,8]$. We train DNN of $10000$ epochs. As an comparison, we choose the Residual neural network with two residual blocks, each block with fully connected layer $D=\{100,100\}$ which gives the result shown in \cref{tab-2} which indicates that the performance of different network structure depends on not only wave number and width, but also the activation function to be selected. In general, the residual network is a little better than that of fully connected network given a good group of representative parameters $\mathcal{K}$ and $\Sigma$.   
	%}
\end{example}

\begin{example}[The 1D cubic Schr\"odinger equation]\label{eq-eg1}
	Here we test the method on the cubic Schr\"odinger equation \cref{Schrodinger-1} in 1D. We use solutions from {the following four} initial conditions to test the accuracy,
	\begin{align}
	%u_0^{\i}(x)&=\exp(-x^2 + i5x),\label{initial-1}\\
	u_0^{\i}(x)&=\left\{\begin{aligned}
	&-0.5 x+1,\quad 0\le x \le 2,\\
	&0.5 x+1,\quad -2\le x < 0,\\
	&0,\qquad \mbox{otherwise},
	\end{aligned}\right. ,\label{initial-1}\\
	u_0^{\ii}(x)&=\sech(x)\exp(i5x),\label{initial-2}\\
	u_0^{\iii}(x)&=\left\{\begin{aligned}
	&\exp(-x^2 + i5x),\quad x \in [-2,2],\\
	&0,\qquad \mbox{otherwise},
	\end{aligned}\right. \label{initial-3}\\
	u_0^{\iv}(x)&=\left\{\begin{aligned}
	&1\quad x \in [-2,2],\\
	&0,\qquad \mbox{otherwise}.
	\end{aligned}\right. \label{initial-4}
	\end{align} 
	
	{The first initial condition is the standard hat function. The last initial condition is a square signal with discontinuities at $x=\pm 2.$ They have no resemblance with the Gaussian wave packets in the training set. }

	As demonstrated in the previous section, to generate data,  the Strang splitting method, combined with the spectral method \cite{bao2002time} are used in the domain $[-\pi \Omega_0,\pi \Omega_0]$ with $\Omega_0=16$ and $N_x$ being the number of Fourier modes. The numerical solution can be captured up to a single-time $T=1$ with $N_t=1000$. The training samples are generated by taking \(\mathcal{K}=\{1,2,\cdots,10\} \textrm{ and } \Sigma=E \) with $h=0.5$. 
	
	We take FCNNs with $D=5$, $m_1=N_x$, $m_2=N_x$ and $M=100$, and train the network for $20000$ epochs. We pick $N_x=8192$ both for the training and testing. The results for solutions from initial conditions {\cref{initial-1} with weak singularity, \cref{initial-2} with a smooth profile, \cref{initial-3} and \cref{initial-4} with discontinuities} are presented in \cref{fig-31}. The relative errors for the density, the real part and the imaginary part of the wave function are found to be  {$6.174\times 10^{-2}$, $1.112\times 10^{-1}$, $9.267\times 10^{-2}$ for \cref{initial-1},  and $4.455\times 10^{-2}$, $5.549\times 10^{-2}$, $7.193\times 10^{-2}$, for \cref{initial-2}, and $1.817\times 10^{-2}$, $2.231\times 10^{-2}$, $2.112\times 10^{-2}$, for \cref{initial-3} and $2.325\times 10^{-1}$, $4.339\times 10^{-1}$, $3.923\times 10^{-1}$, for \cref{initial-4}, respectively, indicating that the smoothness of the solution has an effect on the accuracy of the learned neural network. The last of representations of such solutions in the training set may also be responsible for this outcome. }
	
	%Next, we choose the initial condition as,
	%\begin{align}
	%\psi_0^{\textrm{re}}(x)=sech(x)\cos(5x), \quad 	\psi_0^{\textrm{im}}(x) = sech(x)\sin(5x),
	%\end{align}
	%with $sech(x)=2/(\exp(x)+\exp(-x))$.
	%The results for inputing initial condition \cref{initial-2} are presented in \cref{fig-33}. We then have the relative errors for the density, the real part and the imaginary part of wave function are $6.710\times 10^{-2}$, $1.083\times 10^{-1}$, $1.105\times 10^{-1}$, respectively. The dynamics of wave packets solution with inputing initial condition \cref{initial-2} is presented in \cref{fig-34} of training $2000$ epochs.
	
	\begin{figure}[htbp]
		\centering
		\subfloat{\includegraphics[width=3in]{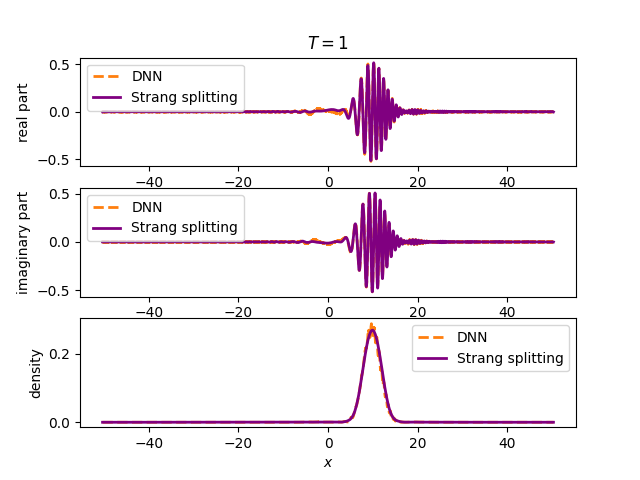}}
		%	\subfloat{\includegraphics[width=3in]{cubic_Schr_error_v2_v1_V2.png}}
		%	\hspace{0.1in}
		\subfloat{\includegraphics[width=3in]{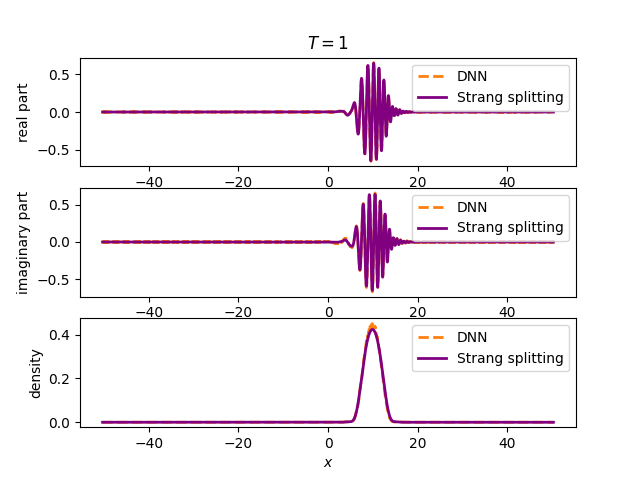}}
		%	\subfloat{\includegraphics[width=3in]{cubic_Schr_error_v4_v1_V1.png}}
		\hspace{0.1in}
		\subfloat{\includegraphics[width=3in]{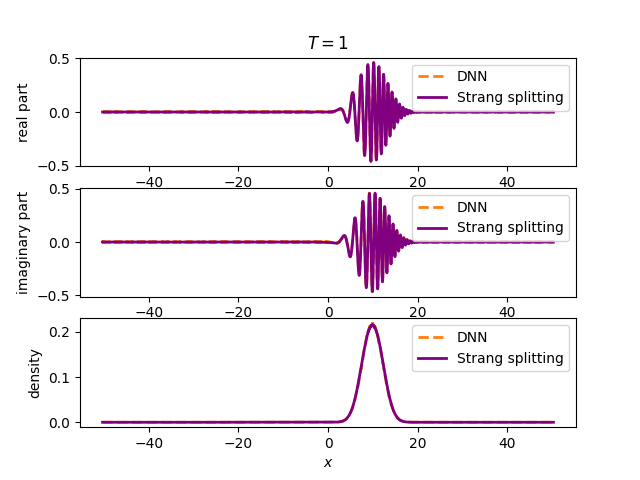}}
		\subfloat{\includegraphics[width=3in]{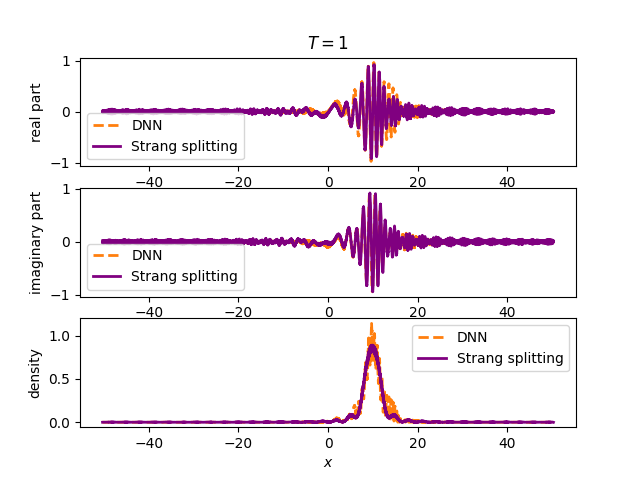}}
		\caption{The learned solution for the 1D cubic NLS \cref{Schrodinger-1}, the direct numerical solution, and error profile at single-time $T=1$ with the initial condition \cref{initial-1} (top left row) and \cref{initial-2} (top right row), \cref{initial-3} (bottom left row) and \cref{initial-4} (bottom right row). The training samples are created with $\mathcal{K}=\{1,2,\cdots,10\}$ and $\Sigma=E$ with $h=0.5$. {Here we take  $\text{tanh}$ as the activation function, and we choose $\lambda=1$. The network is trained for $20000$ iterations.} Left panel: the comparison of the solutions; Right panel: the error.}\label{fig-31}
	\end{figure}

\end{example}

\begin{example}[1D nonlinear Schr\"odinger equation with data from a network]
	
	In most cases, one can not find the exact solution or numerical solution to train the network. We here consider 1D nonlinear Schr\"odinger equation and generate data from a network. This would be useful in cases where analytical or standard numerical methods are difficult to implement.  As a proof-of-concept example, we use the PINNs for \eqref{Schrodinger-1} combined with the CNFD \eqref{eq-CNFD} in time with $f(\rho)=-\rho^{\mu}$. The real part, imaginary part of wave function can be parameterized by a neural network with total residual loss $\mathcal{I}=\mathcal{I}_1+\mathcal{I}_2$ in \eqref{eq-I}. The training samples are generated by PINNs with given initial conditions parametrized by $\mathcal{K}=\{1,2,3,4,5\}$ and various of $\Sigma$. We take PINNs with $D=5$, $m_1=1$, $m_2=2$ and $M=100$, and train the PINNs for $10000$ epochs. We pick $N_x=50$ number of points uniformly sampled in $[-\pi \Omega_0,\pi \Omega_0]$ with $\Omega_0=1$ both for the training and testing purpose. Once the dataset is ready, we take the FCNNs with $D=5$, $m_1=N_x$, $m_2=N_x$ and $M=100$ to train the network for $20000$ epochs.
	{\red For the initial condition \cref{initial-Schr-1},  \cref{fig-nonLinear-onestep_mu} shows the approximation error, for different values of the model parameter $\mu$. An interesting observation, based on the numerical tests  for many activation functions and various choices of $\mathcal{K}$ and $\Sigma$, is that the nonlinearity has an appreciable impact on the accuracy.  } 
	% The stable parameters for activations $\text{relu}$, $\text{tanh}$, $\text{sigmoid}$ and $\text{elu}$ are $E|_{h=0.1}$ with ResNet, $L$ with ResNet, $L$ with ResNet and  $E|_{h=0.5}$, respectively, while sensitive parameters are $L$ with FCNN, $E|_{h=1}$ with ResNet, $E|_{h=0.1}$ with FCNN and $E|_{h=1}$, respectively. As for a fixed group of parameters from training set, the best activations for $E|_{h=0.5}$, $E|_{h=1}$ and $L$ are all $\text{sigmoid}$ with ResNet, while for $E|_{h=0.1}$, taken $\text{relu}$ with ResNet. The sensitive activations for $E|_{h=0.1}$, $E|_{h=0.5}$, $E|_{h=1}$ and $L$ are $\text{sigmoid}$ with FCNN, $\text{relu}$ with FCNN, $\text{tanh}$ with ResNet and $\text{relu}$ with FCNN, respectively. Therefore, the complexity of nonlinear model brings essential difficulty.       
	\begin{figure}[htbp]
		\centering
		\subfloat{\includegraphics[width=3.1in]{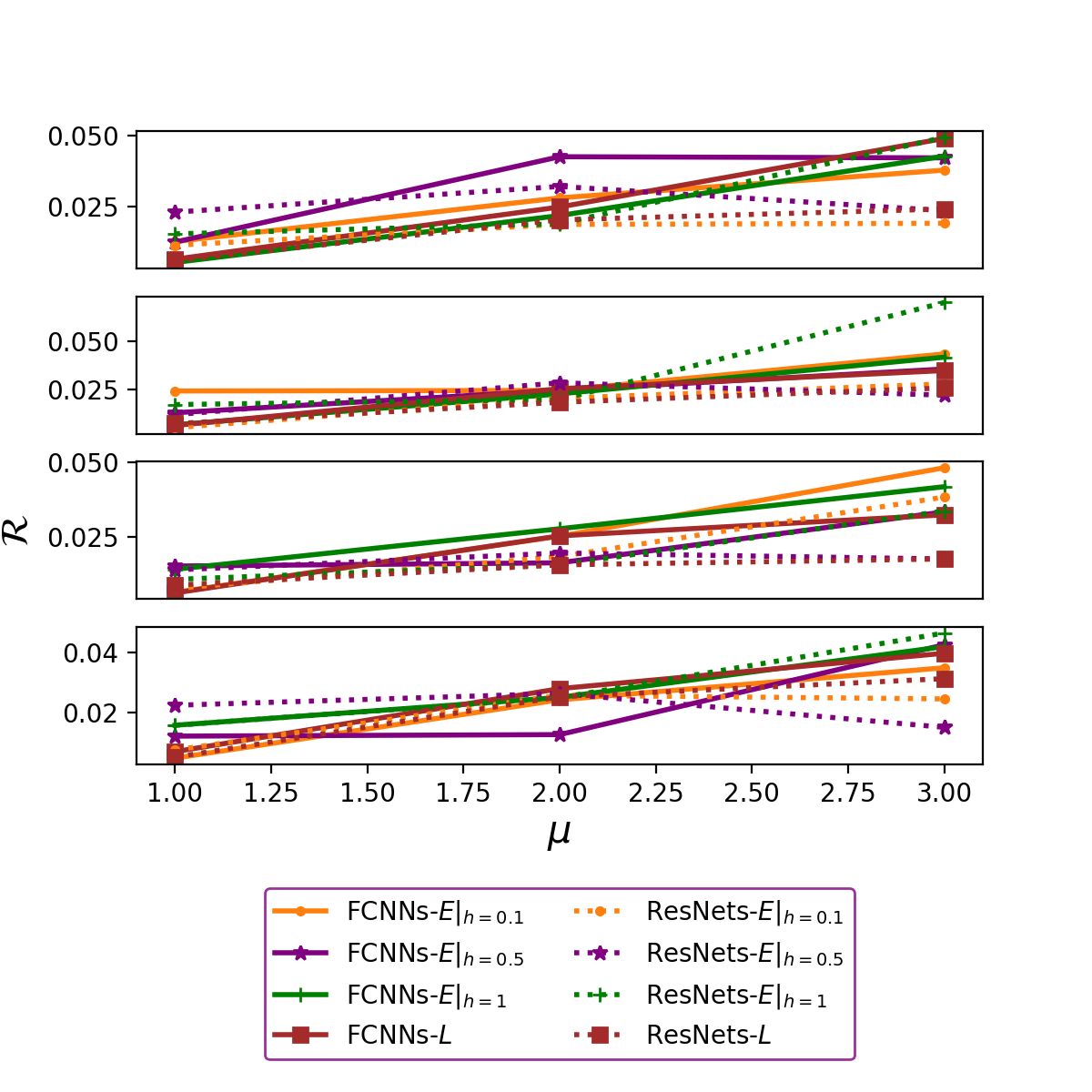}}
		\subfloat{\includegraphics[width=3.1in]{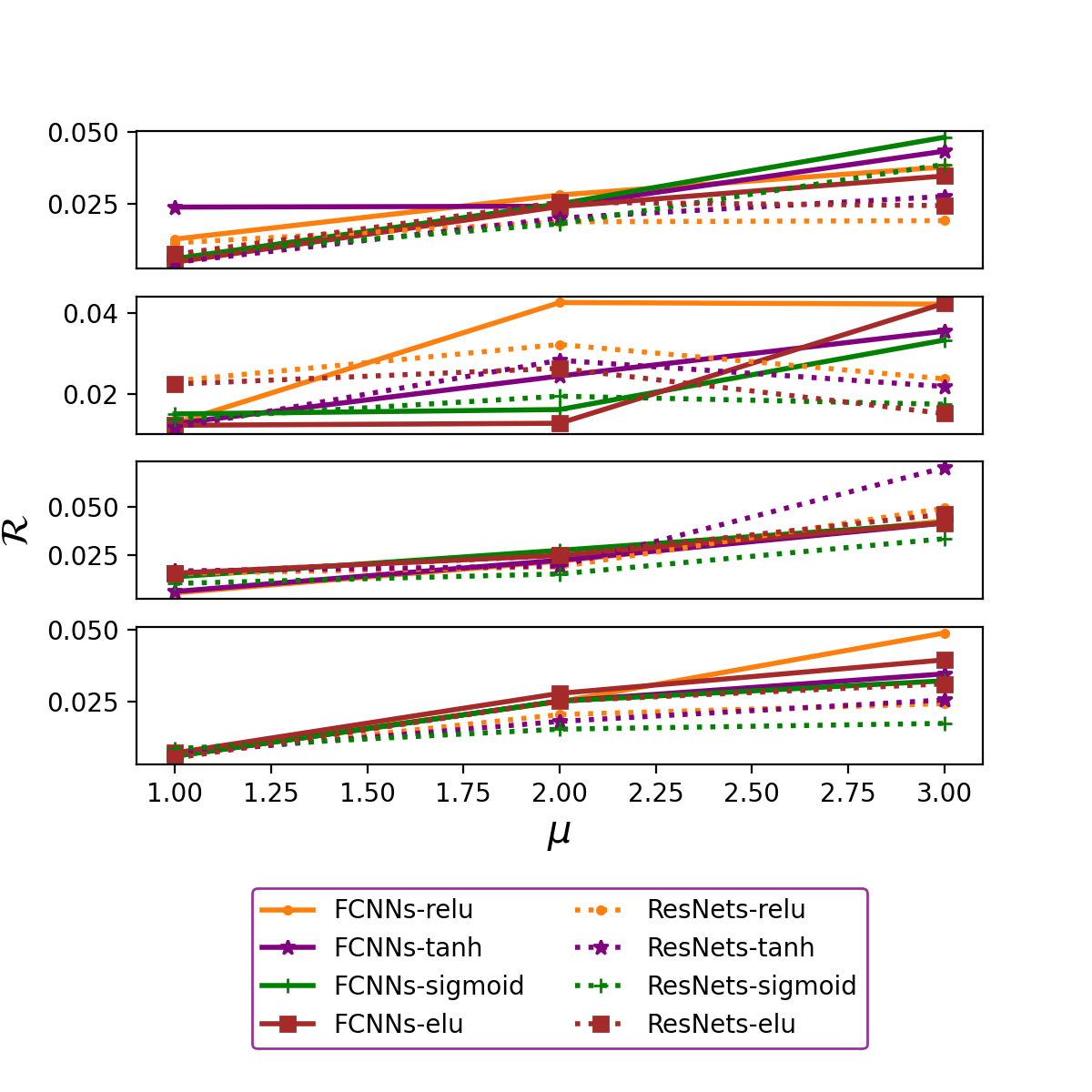}}
		\caption{Approximation error at $T=5\times 10^{-4}$ of the density for the 1D nonlinear Schr\"odinger equation with  initial condition \cref{initial-Schr-1}  using FCNNs and ResNets. Different choices of model parameter $\mu=1,2,3$ and various choices of $\Sigma$. The wave numbers are drawn from $\mathcal{K}=\{1,2,\cdots,5\}$ to build the training set. {Here we take $\lambda=1/4$ and the number of iterations to be $20000$.} Left (from top to bottom): relu, tanh, sigmoid, elu, as the activation functions; Right (from top to bottom): $E|_{h=0.1}$, $E|_{h=0.5}$, $E|_{h=1}$, $L$.}\label{fig-nonLinear-onestep_mu}
	\end{figure}
\end{example}

\begin{example}[1D nonlinear Schr\"odinger equation with a time-dependent potential] 
	\quad 
	
	For problems where physical processes are initiated by an external potential, such as Gross--Pitaevskii equation in Bose--Einstein condensate, {there are two interesting scenarios: (a) mapping the initial condition to the solution at later time given a potential; (b) mapping the potential to the solution at later time given an initial condition.}
	%{\color{blue} It does make sense for learning from potential to solution. But here, in this example, what I actually do is for fixed potential, the mapping is still from initial condition parameterized by $k$ and $\sigma^2$ to solution at later times. However, actually, how to parameterized potential for fixed initial condition, still related to $k$ and $\sigma^2$, right?}
	Motivated by (a), we consider \cref{Schrodinger-1} given a potential that is compactly supported in space. Absorbing boundary condition for this type of problems can be derived \cite{antoine2007review}.
	
	As a specific example, we consider $V(x,t) = E(t)U(x)$,  with time-dependent modulation,
	\begin{align}
	E(t)=E_0\exp(-\gamma(t-t_0)^2)\cos(\omega t). \label{eq-vxt}
	\end{align}
	The spatial part given by,
	\begin{equation}
	U(x)=\begin{cases}
	E_0 x(1-x), \quad x \in [0,1],\\
	0,\quad \textrm{otherwise},
	\end{cases} \label{eq-v}
	\end{equation}
	where $E_0$ denotes the constant intensity of current. In our test, we take $\gamma=1$, $E_0=1$, $t_0=0$ and $\omega=1$. The network setup is the same as that of \cref{eq-eg1} except that we use multiple-time output. As comparison,  the reference solution is computed numerically using the Strang-splitting method combined with a spectral method with the initial condition \cref{initial-2}. The training samples are generated with $\mathcal{K}=\{1,2,\cdots,10\}$ and $\Sigma=E$ with $h=0.5$. The solution learned by the neural network, along with the solution directly computed  are shown in \cref{fig-pot-1}.   The external potential widens the initial wave packets, which subsequently propagate toward the right boundary. The solution represented by the neural network shows great agreement with the direct solution. 
	
	\begin{figure}[ht]
		\centering
		\subfloat[Learned solution with $U(x)$]{\includegraphics[width=3in]{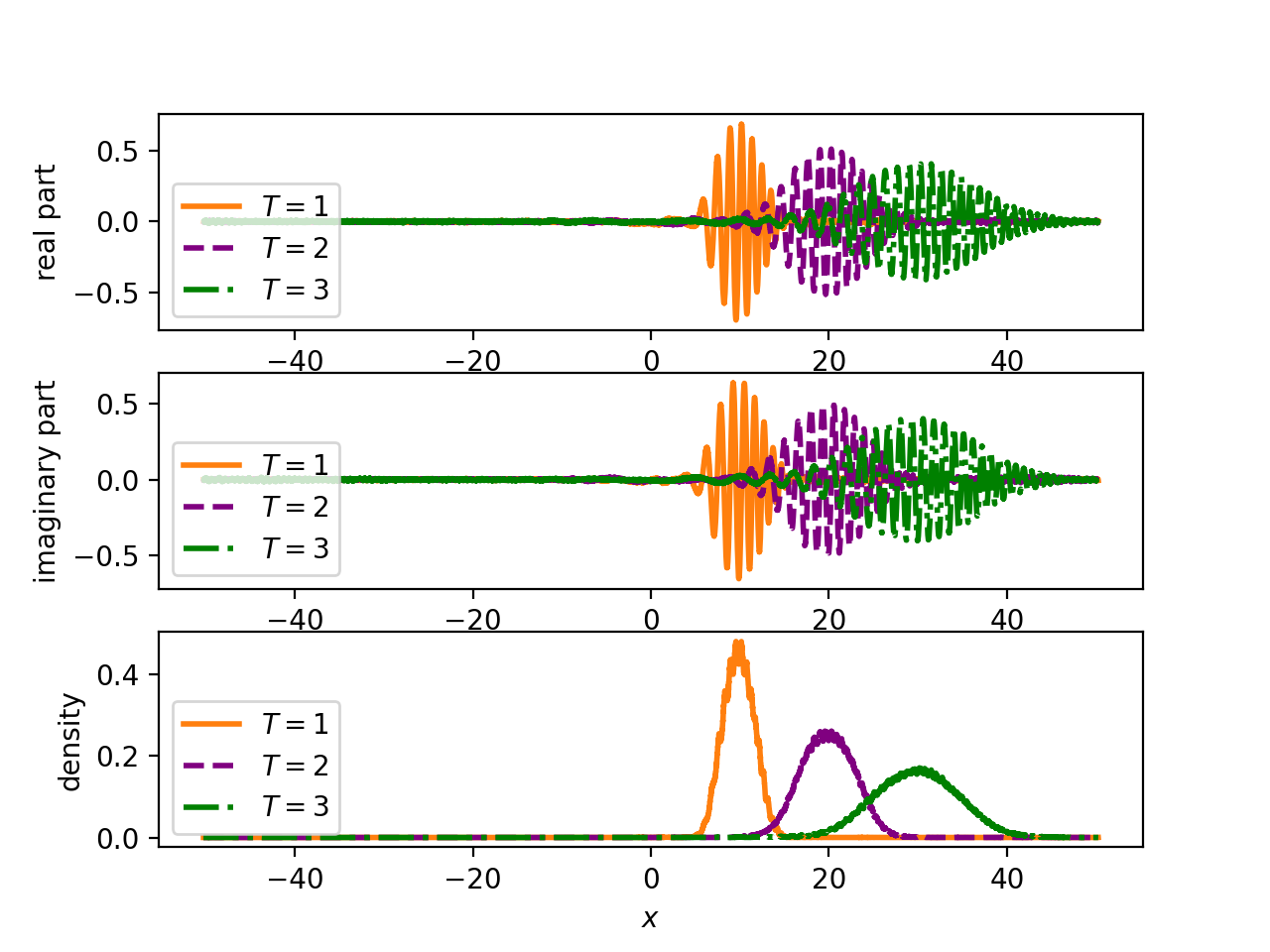}}
		\subfloat[Numerical Strang splitting with $U(x)$]{\includegraphics[width=3in]{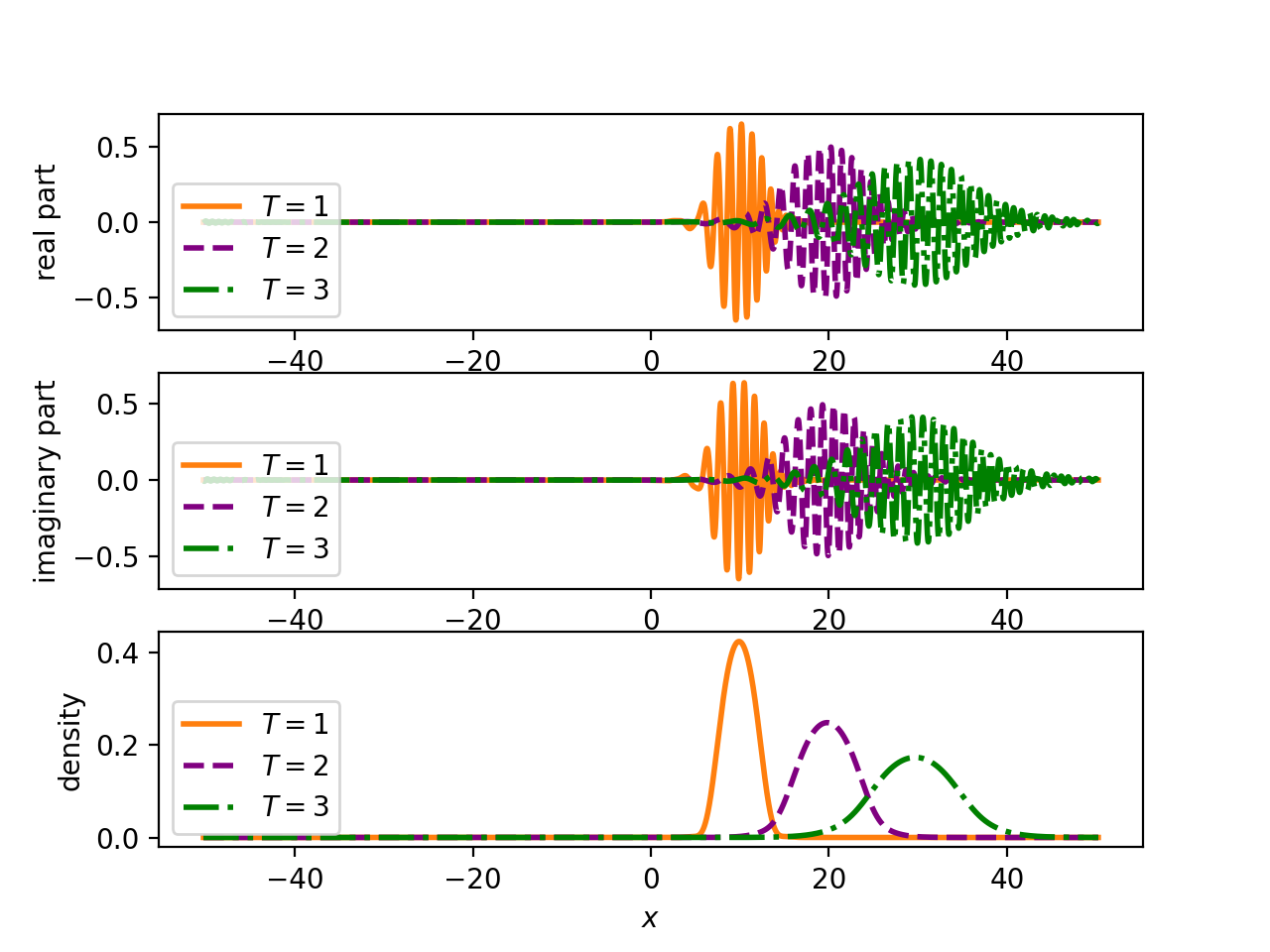}}
		\hspace{0.1in}
		\subfloat[Learned solution with $V(x,t)$]{\includegraphics[width=3in]{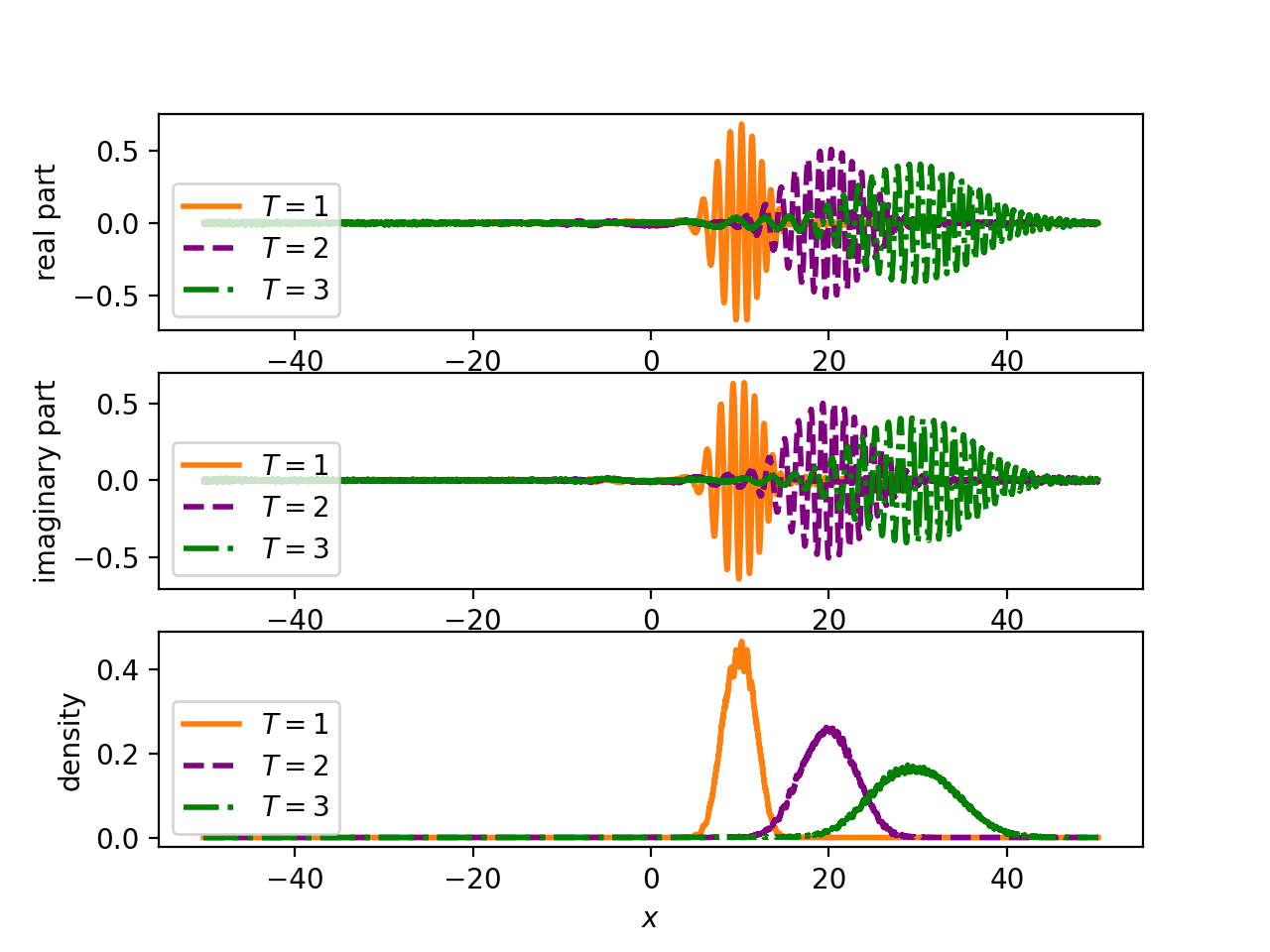}}
		\subfloat[Numerical Strang splitting with $V(x,t)$]{\includegraphics[width=3in]{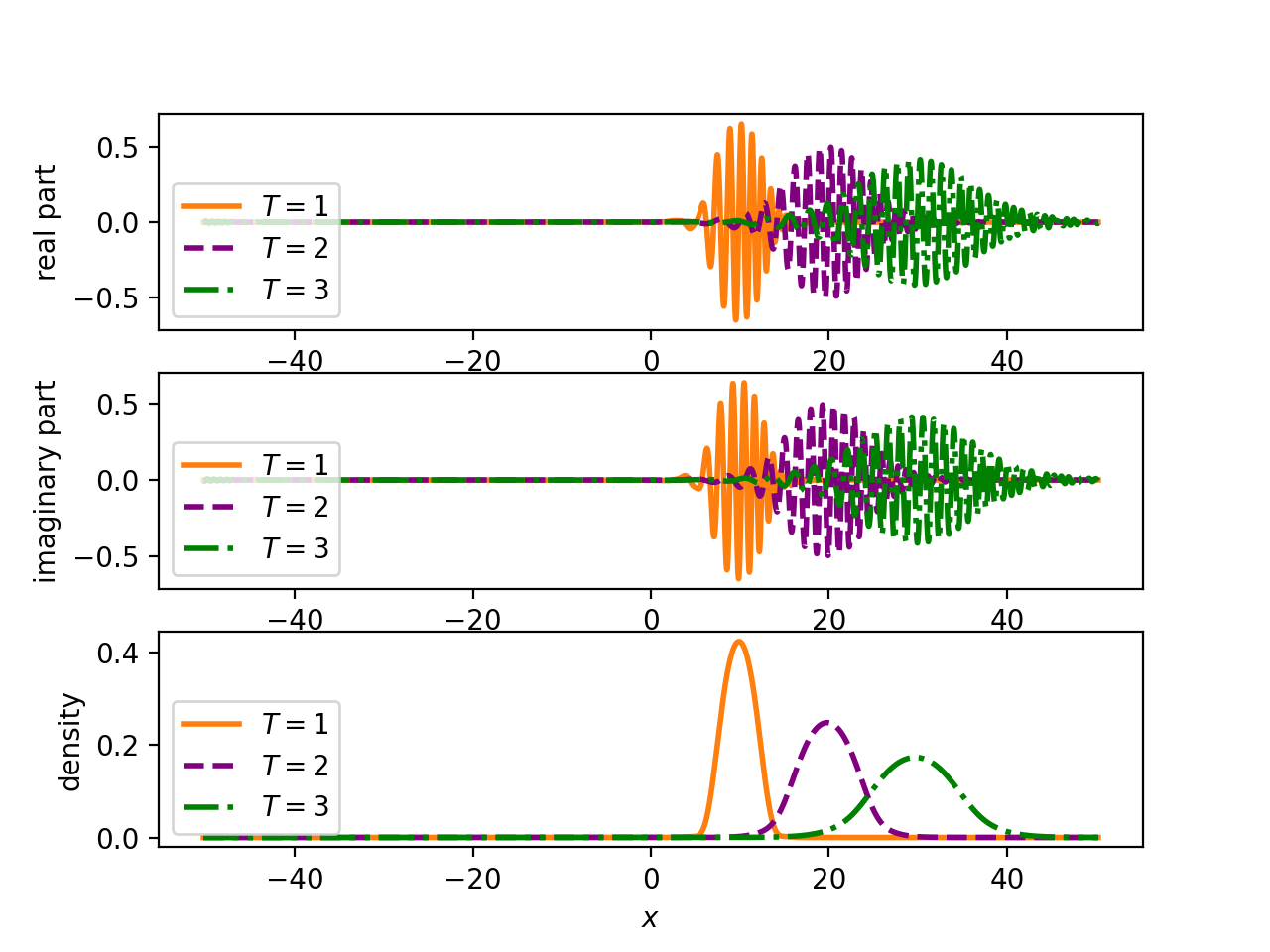}}
		\caption{The learned solution at multiple-time $T=1,2,3$ for the 1D NLS \cref{Schrodinger-1} under an external potential function $V$.   Learned solution (left panel); Numerical solution (right panel). NLS with potential $U(x)$ specified by \cref{eq-v} (top row), and a time-dependent potential $V(x,t)$ specified by \cref{eq-vxt} (bottom row).}\label{fig-pot-1}
	\end{figure}
	
	%\textcolor{blue}{Motivated by (b), we consider \cref{Schrodinger-1} given an initial condition \cref{initial-2} with a compact support. One can choose some %potentials parameterized by the wave number and the band width of Gaussian wave packets as the exactly same as parameterized the initial condition. One %can imagine that the neural network can predict solutions at later time given a newly potential since the mapping learned from potentials to solutions. It is %exactly the case (a) under the time independent potential.}
	
\end{example}

%To study the effectiveness and robustness of this approach for 2D Schr\"odinger equation, the tests here end with 2D nonlinear cubic Schr\"odinger equation without potential. One can generalize the results with other initial conditions endowed with compact support. Same idea can be extended into the higher dimensions for other wave type equations.  

\begin{example}[The 2D cubic Schr\"odinger equation]
	We consider here the  cubic Schr\"odinger equation \cref{Schrodinger-1} in  2D. 
	%The nonlinear Schr\"odinger equation (NLS) in \cite{SULEM1997}.
	%\begin{align}
	%\begin{aligned}\label{eq-1-5}
	%&i\psi_t +\Delta \psi +\alpha|\psi|^2 \psi=0,\quad \x \in \mathbb{R}^2\\
	%&\psi(\x,0)=\psi_0(\x)
	%\end{aligned}
	%\end{align}
	%with cubic nonlinearity and $\alpha=\pm 1$ is an vital paradigm for the envelop dynamics of a weakly nonlinear wave-packet propagating in a dispersive medium. We assume that a positive definite dispersion tensor and concentrate on the case $\alpha=1$ (focusing). 
	The problem is set up as follows: We consider the solution of \cref{Schrodinger-1} in a compact domain $\Omega=[-\Omega_0\pi,\Omega_0\pi]\times [-\Omega_0\pi,\Omega_0\pi]$ with $\Omega_0=2.$ The solutions in $\Omega$ are represented at grid points with $N_x=N_y=64$. The solution in time up to $T=1$ is represented at equally spaced time steps with $N_t=100$. The training and testing are both handled within the domain $\Omega$. For the parameters in the wave packet, we choose $\mathcal{K}= \{1,2,\cdots,5\}$ and we pick $(k_1,k_2) \in \mathcal{K}\times \mathcal{K}$ and 
	\( \Sigma=\{h, 2h, 3h, 4h, 5h, 6h\} \) with $h=0.25$.

	To generate a reference solution, we choose a larger domain $[-2\Omega_0\pi,2\Omega_0\pi]\times [-2\Omega_0\pi,2\Omega_0\pi]$  so that it represents solutions over the entire space within the time period under consideration.  %For the grid points, we choose $N_x=N_y=128$ for the reference solution up to $T=1$ with $N_t=100$ from numerical Strang splitting method and $N_x=N_y=64$ for the training and testing steps with solutions in the subdomain $ [-L\pi/2,L\pi/2]$. 
	For the network, we choose an FCNN with $D=4$, $m_1=N_x N_y$, $m_2=N_x N_y$ and $M=100$. The network is trained for { $40000$} epochs. Then we use the network to predict the solution of the 2D cubic Schr\"odinger equation \cref{Schrodinger-1} with initial condition: 
	{
		\begin{align}\label{eq-5-1}
		\begin{aligned}
		u_0(x_1,x_2)&=\exp[-(x_1^2+x_2^2)+i(3x_1+3x_2)].
		\end{aligned}
		\end{align}
	}
	The results are presented in \cref{fig-29}. The corresponding relative errors for the density, the real  and  imaginary parts of the wave function are { $7.91\times10^{-3}$, $8.98\times10^{-3}$,  and $8.14\times 10^{-3}$}, respectively. The results will improve if we choose $\mathcal{K}=\{1,2,\cdots,5\}$ but $\Sigma=\{h, 1.1h, 1.2h, 1.3h, 1.4h, 1.5h\}$ with $h=0.8$ which give the relative error { $6.48\times 10^{-3}$, $6.53\times 10^{-3}$ and $5.73\times 10^{-3}$} for the density, the real part, and the imaginary parts, respectively. Notice that at time $T=1$, part of the wave packets have moved out of $\Omega$. Therefore, the network has shown an ``absorbing'' property.
	
\end{example}

\begin{example}[Wave propagation in irregular domains]
	
	{ We extend the previous example to irregular domains. Specifically, 
		we consider a circular disk, and an $L$-shape domain. This does not impose further difficulty on the method: We simply repeat the training procedure over the corresponding domains. As shown in \cref{fig-30} for the disk domain, the corresponding relative errors for the density, the real and imaginary parts of the wave function are $1.21\times 10^{-2}$, $9.99\times 10^{-3}$, $8.28\times 10^{-3}$ for $\mathcal{K}=\{1,2,3,4,5\}$, $\Sigma=\{h,2h,3h,4h,5h,6h\}$ with $h=0.25$ and $4.88\times 10^{-3}$, $8.08\times 10^{-3}$, $8.01\times 10^{-3}$ for $\Sigma=\{h,1.1h,1.2h,1.3h,1.4h,1.5h\}$ with $h=0.8$, respectively. As shown in \cref{fig-32} for the L-shape domain, the corresponding relative errors for the density, the real and imaginary parts of the wave function are $7.86\times 10^{-3}$, $9.45\times 10^{-3}$, $8.77\times 10^{-3}$ for $\mathcal{K}=\{1,2,3,4,5\}$, $\Sigma=\{h,2h,3h,4h,5h,6h\}$ with $h=0.25$ and $9.46\times 10^{-3}$, $1.18\times 10^{-2}$, $8.94\times 10^{-3}$ for $\Sigma=\{h,1.1h,1.2h,1.3h,1.4h,1.5h\}$ with $h=0.8$, respectively. } In both cases, the neural network has demonstrated an absorbing property, allowing the waves to propagate out of a  domain with complex geometry.  
	
	\begin{figure}[htbp]
		\centering
		\subfloat{\includegraphics[width=2in]{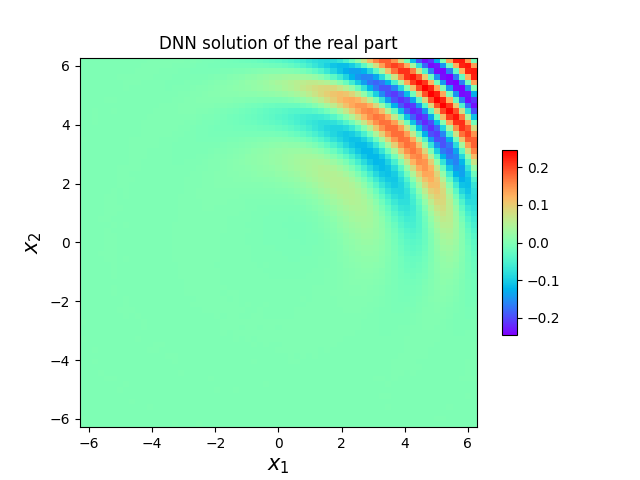}}
		\subfloat{\includegraphics[width=2in]{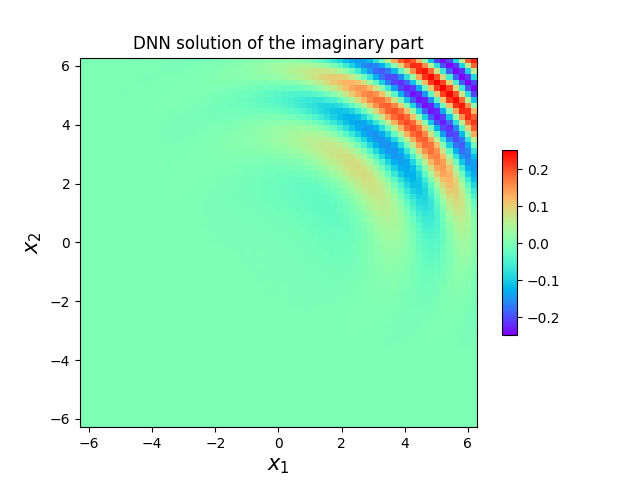}}
		\subfloat{\includegraphics[width=2in]{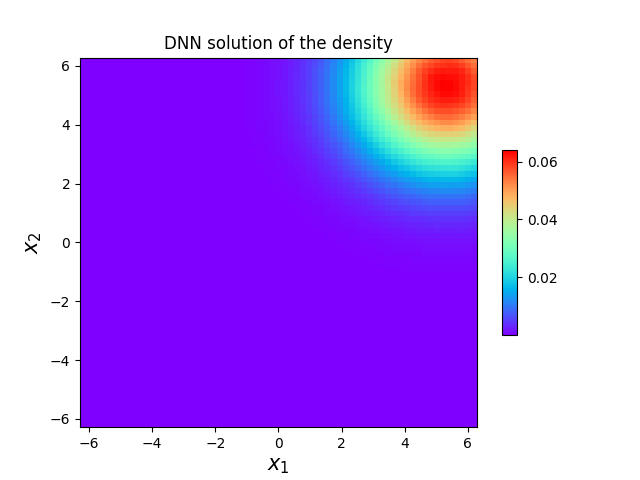}}
		\hspace{0.1in}
		\subfloat{\includegraphics[width=2in]{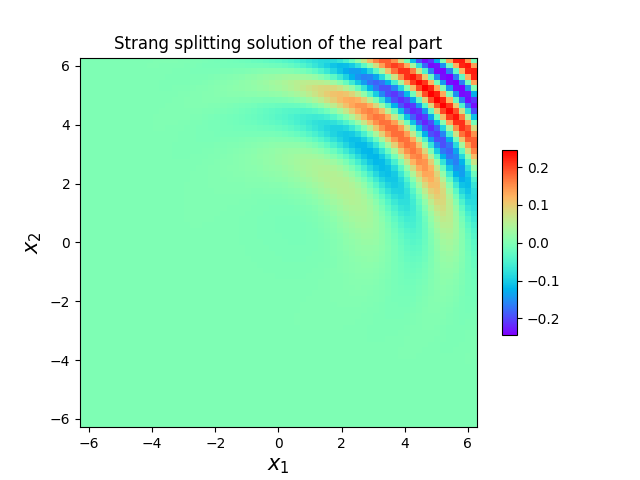}}
		%	\subfloat{\includegraphics[width=2in]{error_re_t_1_v5_v1_v1_V1.png}}
		\subfloat{\includegraphics[width=2in]{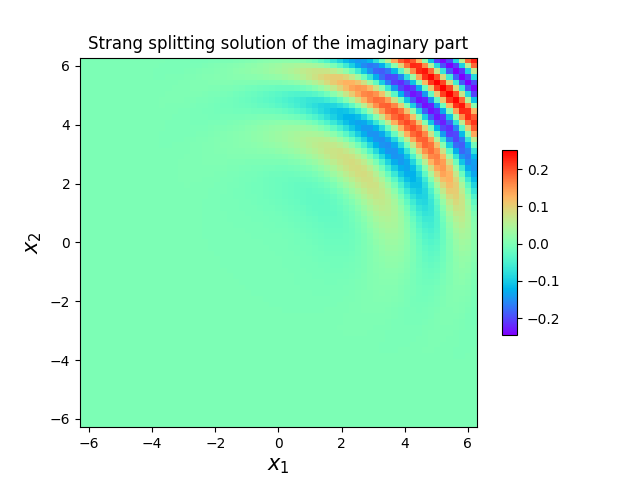}}
		%	\subfloat{\includegraphics[width=2in]{error_im_t_1_v5_v1_v1_V1.png}}
		\subfloat{\includegraphics[width=2in]{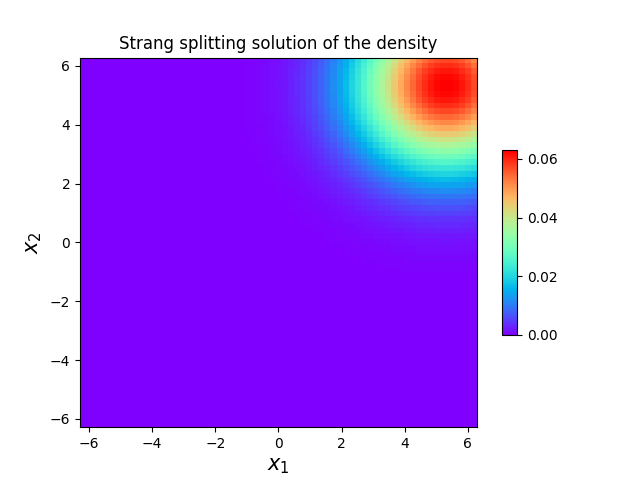}}
		%	\subfloat{\includegraphics[width=2in]{error_rho_t_1_v5_v1_v1_V1.png}}
		\caption{Solution of the 2D cubic Schr\"odinger equation in a square domain by the DNN (Top row), direct numerical method (Bottom row), evaluated at time $T=1$ for the initial condition \cref{eq-5-1}. {A network with  activation function $\text{tanh}$ is trained for $40000$ iterations to represent the solutions.} Left panel: the real part of wave packet solution; Middle panel: the imaginary part; Right panel:  the electron density $\rho$.}\label{fig-29}
	\end{figure}
	
	\begin{figure}[htbp]
		\centering
		\subfloat{\includegraphics[width=2.1in]{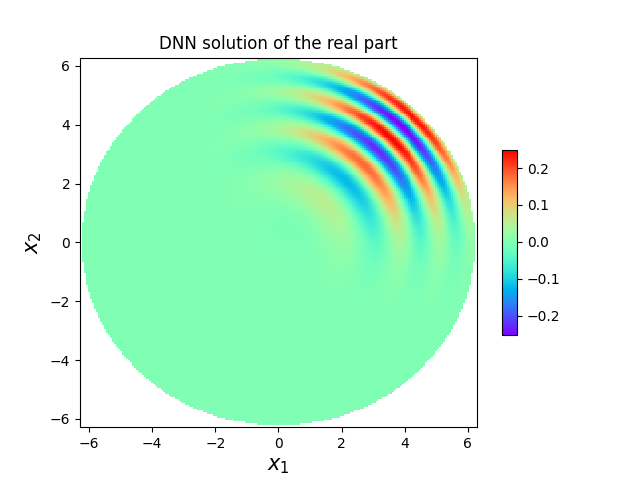}}
		\subfloat{\includegraphics[width=2.1in]{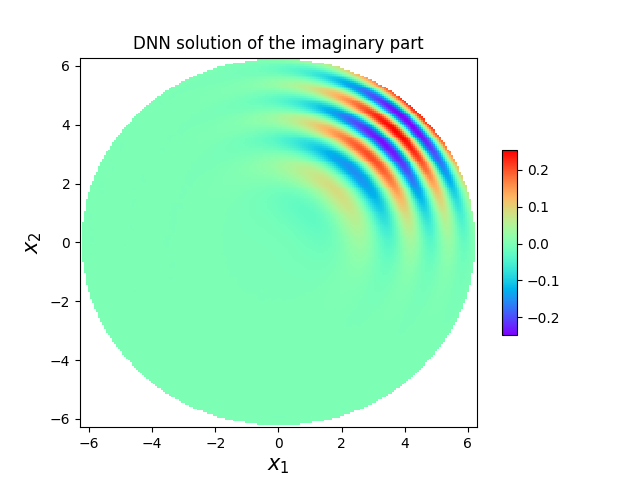}}
		\subfloat{\includegraphics[width=2.1in]{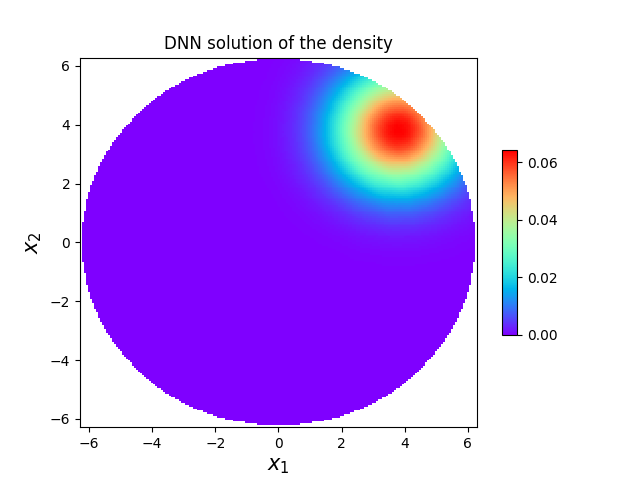}}
		\hspace{0.1in}
		\subfloat{\includegraphics[width=2.1in]{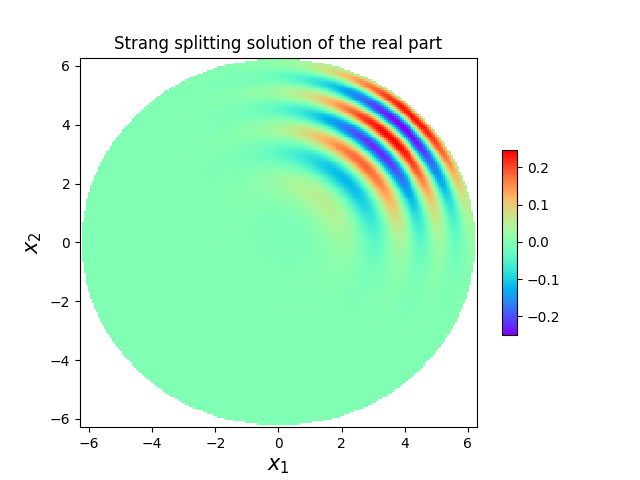}}
		%	\subfloat{\includegraphics[width=2in]{error_re_t_1_v5_v1_v1_V1.png}}
		\subfloat{\includegraphics[width=2.1in]{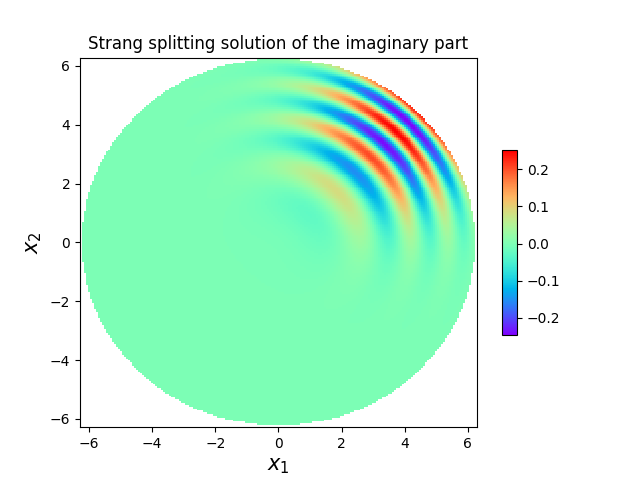}}
		%	\subfloat{\includegraphics[width=2in]{error_im_t_1_v5_v1_v1_V1.png}}
		\subfloat{\includegraphics[width=2.1in]{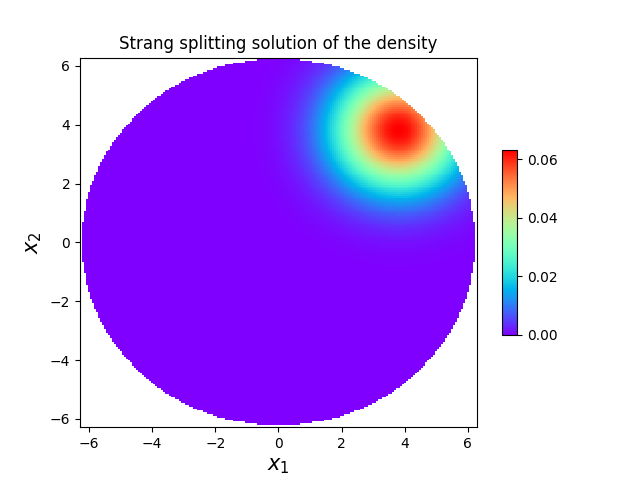}}
		%	\subfloat{\includegraphics[width=2in]{error_rho_t_1_v5_v1_v1_V1.png}}
		\caption{Solution of the 2D cubic Schr\"odinger equation in a circular domain by the DNN (Top row), direct numerical method (Bottom row), evaluated at time $T=1$ for the initial condition \cref{eq-5-1}. {Here the activation function $\text{tanh}$ is used and the training takes $40000$ iterations.} Left panel: the real part of wave packet solution; Middle panel: the imaginary part; Right panel:  the electron density $\rho$.}\label{fig-30}
	\end{figure}
	
	\begin{figure}[htbp]
		\centering
		\subfloat{\includegraphics[width=2.1in]{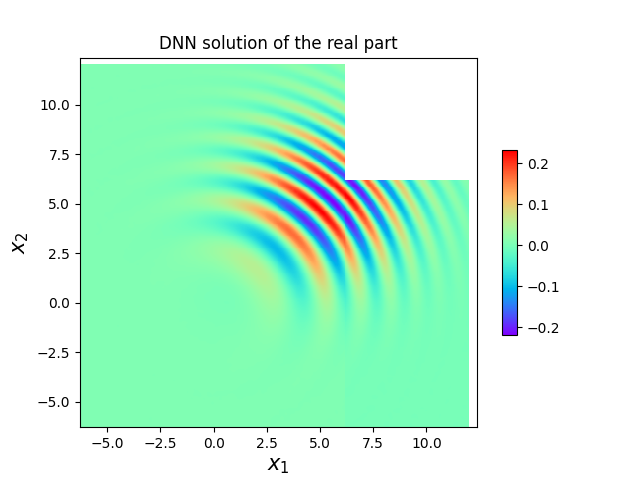}}
		\subfloat{\includegraphics[width=2.1in]{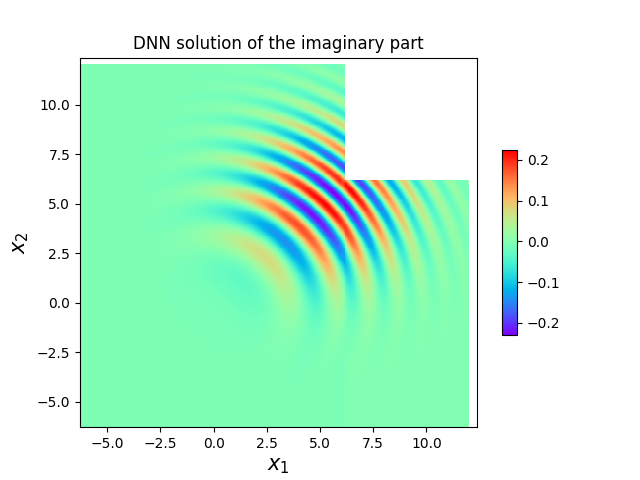}}
		\subfloat{\includegraphics[width=2.1in]{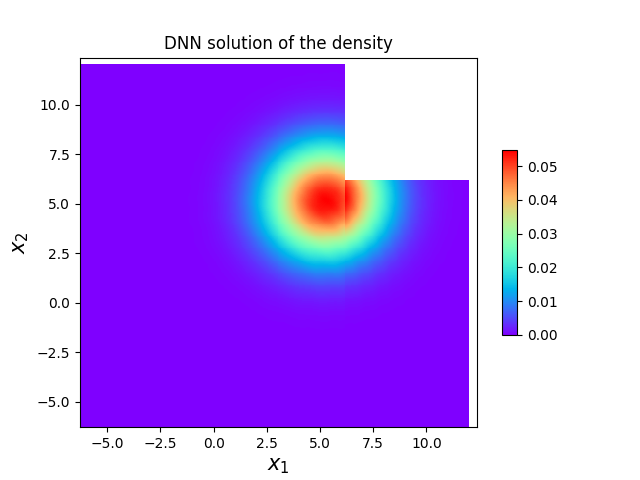}}
		\hspace{0.1in}
		\subfloat{\includegraphics[width=2.1in]{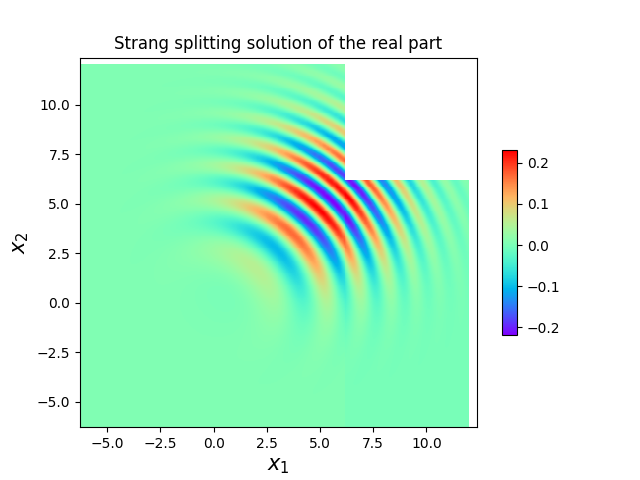}}
		%	\subfloat{\includegraphics[width=2in]{error_re_t_1_v5_v1_v1_V1.png}}
		\subfloat{\includegraphics[width=2.1in]{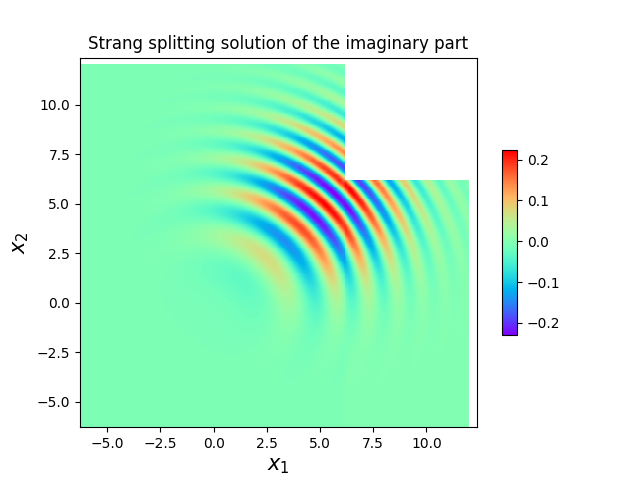}}
		%	\subfloat{\includegraphics[width=2in]{error_im_t_1_v5_v1_v1_V1.png}}
		\subfloat{\includegraphics[width=2.1in]{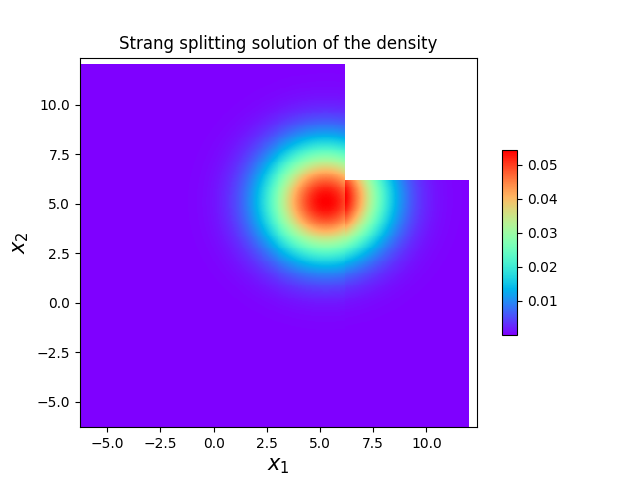}}
		%	\subfloat{\includegraphics[width=2in]{error_rho_t_1_v5_v1_v1_V1.png}}
		\caption{Solution of the 2D cubic Schr\"odinger equation in an L-shape  domain by the DNN (Top row), direct numerical method (Bottom row), evaluated at time $T=1$ for the initial condition \cref{eq-5-1}. {Here the activation function $\text{tanh}$ is used and the network is trained for  $40000$ iterations.} Left panel: the real part of wave packet solution; Middle panel: the imaginary part; Right panel:  the electron density $\rho$.}\label{fig-32}
	\end{figure}
	
\end{example}

\section{Discussion and Conclusions}
\label{sec:conclusions}

In this paper, we have proposed a machine-learning method to solve wave equations over unbounded domains without introducing artificial boundary conditions. As examples, we considered the Schr\"{o}dinger equation and the second-order acoustic wave equation. Results show that the proposed method has good interpolative accuracy and some extrapolative accuracy. 

All simulations are implemented in MacBook Pro Intel Core i5 ($1$ CPU, $4$ Kernels and $8$ Gb random access memory). The method provides an alternative for finite-time simulation of wave propagation. On the other hand, we found that wave propagations over long time period still remains a challenge for the neural network approximation. %This will be of particular interest and shall be explored in future.

Unlike conventional numerical methods for solving wave equations, e.g., \cites{diaz2006time,bao2005convergence}, results for rigorous error bounds from neural network approximations of PDEs over unbounded domains are scarce. Therefore we rely on extensive numerical experiments and we report some direct observations here.

\smallskip

\noindent {\bf Neural network as a PDE solver.}  Recently, the machine leaning approach has been applied to wave equations in \cites{cai2019phase}. It is important to point out  though that most of those effort focused on solving PDEs using  neural networks on {\it bounded} domains \cites{lagaris1998artificial,khoo2017solving,sirignano2018dgm,Lyu2020Jun}, while the current work is aimed at representing solutions of time-dependent hyperbolic  PDEs on unbounded domains. 

\medskip 
\noindent {\bf The choice of the data set and network structures.}  
It is apparent from the numerical tests that the accuracy of the neural network representation depends crucially on the choice of the training set, as well as the network structure and the activation functions.  {When the test set lies in the training sets (or the range)}, we interpret the representation \eqref{eq: nnet0} as an interpolation.  The accuracy is generally satisfactory with the error on the scale of $10^{-3}$. { We observe that the $\tanh$ function for the second-order acoustic wave equation stood out as the best choice with robust overall performance in both FCNN and ResNet, as suggested by the comparison in Table \ref{tab-2-1}, while the $\relu$ function for the Schr\"odinger equation is the best choice in Table \ref{tab-2}.}
%the $\relu$ activation function stood out as the best choice with robust overall performance  in both FCNN and ResNet, as suggested by the comparison in Tables \ref{tab-2-1}, \ref{tab-2-2}, \ref{tab-2}, with the exception of the results in \ref{tab-3}. 
The comparison between the performance of FCNN and ResNet seems more subtle; { The ResNet seems to perform better in more cases.}
%seems to depend on the choice of the activation function and the specific test problem.
% in most of these cases, with the exception that for the  $\relu$ activation function, ResNet offers slightly better accuracy than FCNN with exponential spacing for $\Sigma$. When $\relu$ is used, the set $\Sigma$ with linear spacing is slightly better than a set with exponential spacing. 

In contrast to interpolations, extrapolations can arise in different ways. For instance, the initial condition can be of a different function form than those in the training set. Results in \cref{fig-31} are such examples.  FCNNs seem to offer consistent results in this case. An  extrapolation also comes up when the network \eqref{eq: nnet0} is used to predict solutions at a later time. We observed from \cref{tab-extra-schr-t} that 
the accuracy of such an approximation can be guaranteed for a short time period, and the FCNN with { $\elu$} activation function yields much better results than other choices.   Another scenario of an extrapolation is
when the wave number (or the width of the packet) in the initial condition is outside the range of the set $\mathcal{K}$ (or $\Sigma$) that is associated with the training set. Generally, the error grows when the wave number is further away from the set $\mathcal{K}$. But for a specific case, the results among different choices of the activation functions are mixed. The activation function $\elu$ seems to give reasonable accuracy in all the cases tested. 
%in this case is mixed: For the second order wave equation, under the same pre-selected $\mathcal{K}$ and $\Sigma$, ResNet is better than FCNN for this purpose when the activation function is chosen as $\relu$, $\tanh$ or $\sigmoid$. In the case of selecting $\elu$ as an activation function, FCNN is slightly better than ResNet. For  the Schr\"odinger equation,
% the results are mixed. %  as the learning model, FCNN is slightly better than ResNet for selected activation functions $\relu$ and $\tanh$. In this case, $\relu$ is a good choice for FCNN. The activation function $\elu$ is a good choice for ResNet. 

%\textit{In terms of the performance of interpolation with respect to the wave number out of training sets}, ResNet is better than FCNN under the selecting activation functions $\tanh$ or $\sigmoid$. Specifically for the model of the second order wave equation, ResNet is also slightly better than FCNN under $\relu$ and $\elu$. The activation function $\elu$ is a good choice for this model. However, once the Schr\"odinger equation is considered, $\tanh$ is a good choice. In this case, FCNN is slightly better than ResNet under activation functions $\elu$ or $\relu$. 

%is this aspect is slightly better overall than ResNet.% with the activation functions $\tanh$, $\relu$, and $\elu$. In contrast, ResNet is slightly better than FCNN as $\sigmoid$ selected. Among these activation functions, $\tanh$ yields the best results.    

Due to the wave propagation nature, we proposed to use wave packets to create the training set. So far, our numerical tests have not singled out 
an optimal strategy. Although larger selections of    $\mathcal{K}$ and $\Sigma$  generally give better results, they inevitably lead to larger training dataset. 
One possible direction is to start with a larger set of training data, and then use the proper orthogonal decomposition (POD) to extract the most relevant basis. In high dimensions, we might use the (quasi-) Monte Carlo methods to handle and extract the representative elements of $\Sigma$ and $\mathcal{K}$. 

The current approach  excludes nonlinear waves, e.g., shock waves \cite{dafermos2005hyperbolic}, contact discontinuities, solitons \cite{forest1986geometry}, etc. It would be interesting to investigate the performance of the neural network in those scenarios as well. 

\medskip

\noindent {\bf The relation to absorbing boundary conditions.} The current approach targets the same type of problems as absorbing boundary conditions, viz., wave propagation processes  that occur in an unbounded domain, but are triggered by initial conditions or external signals that are localized in a bounded domain. {However, rather than using the neural network to incorporate the absorbing boundary condition into the FDTD procedure  \cites{chen2020learning,yao2018machine,yao2020enhanced,AAMM-12-1384}, which involves the history of solutions at the boundary,} we directly map the initial condition to the solution at time instances of interest. The numerical results suggest that such an approximation also exhibits absorbing properties. 
%{One can solve the wave propagation over unbounded domain using the truncated higher order ABC, which is an approximately method. One can construct the exact Dirichlet-to-Neumann (DtN) mapping to solve the original unbounded problems by usage of neural network formulation, which shall be explored in the future work.} 
In addition to the wave equations we discussed in this paper, the results suggest that this framework can be extended to other wave propagation problems, e.g., those from fluid mechanics \cite{han1996artificial}, elasticity \cites{becache2001fictitious,guddati2000continued}, and molecular dynamics  \cites{karpov2005green,li2006variational,weinan2001matching}.

\medskip

%\noindent {\bf The relation to semi-classical Schr\"odinger equation.} In the numerical simulations of semi-classical Schr\"odinger equation, one is interested in the solution behaviour as $\epsilon\to 0$. Empirically, the parameter $\epsilon$ smaller, the simulation more difficult. We expect the neural network can relax or overcome this barrier by learning end-to-end. 

\noindent {\bf High-dimensional problems. } One of the distinct advantage of neural networks is the ability to treat high-dimensional problems. This has been demonstrated in various type of PDEs. A potential application of the current approach is to many-particle Schr\"odinger equations. For instance, in ionization problems, electrons can be driven away from nuclei via a laser field, and traditionally, such problems have been treated using effective models and absorbing boundary conditions, e.g., in the context of time-dependent density-functional theory \cite{yabana2006real}. This work is currently underway.

\appendix

\section*{Acknowledgments}
This work is supported in part by the financial support from the program of China Scholarships Council No. 201906920043 (C. Xie), National Science Foundation of China Grant No. 11971021 (J. Chen).

\bibliographystyle{amsplain}
\bibliography{references}
%\begin{thebibliography}{10}
%
%\bibitem {A} T. Aoki, \textit{Calcul exponentiel des op\'erateurs
%microdifferentiels d'ordre infini.} I, Ann. Inst. Fourier (Grenoble)
%\textbf{33} (1983), 227--250.
%
%\bibitem {B} R. Brown, \textit{On a conjecture of Dirichlet},
%Amer. Math. Soc., Providence, RI, 1993.
%
%\bibitem {D} R. A. DeVore, \textit{Approximation of functions},
%Proc. Sympos. Appl. Math., vol. 36,
%Amer. Math. Soc., Providence, RI, 1986, pp. 34--56.
%
%\end{thebibliography}

\end{document}